\documentclass[final,1p,times]{elsarticle}
\usepackage{algorithm}
\usepackage{algorithmic}
\usepackage[algo2e]{algorithm2e}
\usepackage{amssymb}
\usepackage{amsthm}
\usepackage{amscd}
\usepackage{amsmath}
\usepackage{amsfonts}
\usepackage{amssymb}
\usepackage{booktabs}
\usepackage{graphicx}
\usepackage{color}
\usepackage{framed}
\usepackage{comment}
\usepackage{chngcntr}
\numberwithin{equation}{section}
\newtheorem{theorem}{Theorem}[section]

\newtheorem{lemma}[theorem]{Lemma}
\newtheorem{corollary}[theorem]{Corollary}

\usepackage{mathrsfs}

\usepackage{titletoc}
\biboptions{sort&compress}
\usepackage{subfigure}
\usepackage{hyperref}
\usepackage[pagewise,mathlines]{lineno}
\usepackage{geometry}

\usepackage{bm}

\biboptions{sort&compress}

\makeatletter

\newcommand{\Rmnum}[1]{\expandafter\@slowromancap\romannumeral #1@}
\makeatother

\usepackage{bbding}

\journal{***}

\begin{document}

\begin{frontmatter}

\title{A modified Ricci flow on arbitrary weighted graph}
		
\author[author1]{Jicheng Ma}
\ead{2019202433@ruc.edu.cn}
\author[author3]{Yunyan Yang{\footnote{corresponding author}}}
\ead{yunyanyang@ruc.edu.cn}

\address{$^1$School of Mathematics, Renmin University of China, Beijing, 100872, China}	

\begin{abstract}

In this paper, we propose a modified Ricci flow,  as well as a quasi-normalized Ricci flow, on arbitrary weighted graph. Each of these two flows has a unique global solution. In particular, these global existence and uniqueness results do not require an exit condition proposed by
Bai et al in a recent work \cite{Bai-Lin}. As applications, these two Ricci flows are applied to community detection for complex networks, including Karate Club, American football games, Facebook, as well as artificial networks.
In our algorithms, unlike in \cite{Ni-Lin,Lai X}, there is no need to perform surgery at every iteration, only one surgery needs to be performed
after the last iteration. From three commonly used criteria for evaluating community detection algorithms, ARI, NMI and Q, we conclude that our algorithms outperform existing algorithms, including Ollivier's Ricci flow \cite{Ni-Lin}, normalized Ollivier's Ricci flow and normalized Lin-Lu-Yau's Ricci flow \cite{Lai X}. The codes for our algorithms are available at https://github.com/mjc191812/Modified-Ricci-Flow.
\end{abstract}

\begin{keyword}
weighted graph; Ricci curvature; Ricci flow; community detection
\\
\MSC[2020] 05C21; 05C85; 35R02; 68Q06
\end{keyword}
		
\end{frontmatter}
	
\section{Introduction}

Ricci flow was first introduced by Hamilton \cite{Hamilton} in 1982. It was originally designed to deform the Riemannian metric on a smooth Riemannian manifold $(M,g_0)$, and was explicitly expressed as an evolution equation
$$\left
\{\begin{array}{ll}
\partial_t g=-2{\rm Ric}\\[1.2ex]
g(0)=g_0,
\end{array}
\right.
$$
 where ${\rm Ric}$ stands for the Ricci curvature on Riemannian manifold $(M,g(t))$. The power of this flow has already been known to the world. For example, it was used by Brendle-Schoen \cite{Brendle-Schoen} to prove the differential sphere theorem, and also used by Perelman \cite{Perelman} to solve the Poincar\'e conjecture.

Intuitively, Ricci flow is valuable for understanding the evolution and community structure of networks. One can think of a network as a discretization of a high-dimensional manifold, similar to a 3-manifold, and the communities in the network as components in the geometric decomposition of the 3-manifold. Perelman's work \cite{Perelman} has shown that Ricci flow can predict the geometric components of a 3-manifold, suggesting that a discrete Ricci flow on a network should be able to detect the community structure. Analogous to the work of Hamilton and Perelman on Ricci flow, the number of iterations and the threshold value for surgery in the flow process may depend on individual networks.

Indeed, as expected, Ricci flow performs well in complex networks. In 2019, based on Ollivier's Ricci flow  \cite{Ollivier-1,Ollivier-2} and a surgery procedure, Ni et al \cite{Ni-Lin} provided an outstanding community detection method. Later, the same problem was re-solved by Lai et al \cite{Lai X} by using a normalized Ricci flow, which is based on Lin-Lu-Yau's Ricci curvature \cite{Lin-Lu-Yau} and a star coupling Ricci curvature of Bai et al \cite{Bai-Huang}. {{The star coupling Ricci curvature is a computational improvement that avoids limit operations.}} Let's say a few more words about complex networks, which are commonly used to represent connections between elements in various fields, including social networks \cite{Wasserman S}, biochemistry (such as protein-protein \cite{Bhowmick} networks, metabolic networks, and gene networks), and computer science \cite{Ni C C,Tauro S L}. It's widely known that many real-world networks exhibit community structures, where nodes within the same community are closely connected, while nodes from different communities are sparsely connected. Recognizing these community structures is essential for identifying key functional components and supporting processes on networks, such as disease spread, information dissemination, and behavioral patterns. Lots of algorithms \cite{Yang-Algesheimer,Fortunato,Newman M E J,Leskovec,Clauset-Newman-Moore,Peel,Girvan M} have been developed to detect and separate communities. Most of these algorithms focus on identifying dense clusters in a graph, using randomized approaches like label propagation or random walks, optimizing centrality measures such as betweenness centrality, or considering measures like modularity. In contrast to these methods, \cite{Ollivier-1,Ollivier-2,Ni-Lin,Lai X} suggest discrete Ricci flows on weighted graphs, which offers broader applicability and greater solution stability.

More precisely, in \cite{Ollivier-1}, Ollivier suggested using the following equation as Ricci flow with continuous time parameter $t$: for each edge $e\in E$, the weight $w_e(t)$ satisfies the ordinary differential system
\begin{equation}\left\{\begin{array}{lll}\label{flow}
\frac{d}{dt}w_e(t)=-\kappa_e(t)w_e(t)\\[1.2ex]
w_e(0)=w_{0,e},\end{array}
\right.
\end{equation}
where $w_{0,e}$ denotes the initial weight and $\kappa_e(t)$ stands for Ollivier's Ricci curvature of the edge $e$ on a connected
weighted graph $(V,E,w(t))$.   For a discrete Ricci flow, in each iteration, the process generates a time dependent family of weighted graph $(V,E,w(t))$ such
that each weight $w_e(t)$ changes proportional to  $\kappa_e(t)$ at time $t$.  In \cite{Ni-Lin}, slightly different from the discrete version of (\ref{flow}), Ni et al used the iteration system
\begin{equation}\label{Ni00}
w_e^{(i+1)}=\rho_e^{(i)}-s\kappa_e^{(i)}\rho_e^{(i)}
\end{equation}
in their algorithm,
where $w_e^{(i)}$, $\kappa_e^{(i)}$ and $\rho_e^{(i)}$ denote the weight, Ollivier's Ricci curvature and the distance, respectively, of the edge $e$ at the $i$th iteration, and $s$ is a step size. Here, if one represents $e=xy$, then the distance is represented by
$$\rho_e^{(i)}=\inf_{\gamma\in\Gamma_{e}}\sum_{\tau\in\gamma}w_\tau^{(i)},$$
where $\Gamma_{e}$ is a set of all paths connecting vertices $x$ and $y$.
It is conceivable that if the volume of a weighted graph tends to zero along the Ricci flow, then the community detection effect of the network would not be satisfactory. To overcome this shortcoming, Lai et al \cite{Lai X} proposed a normalized discrete Ricci flow, the $i$th process of which says
\begin{equation}\label{Lai-1-0}
w_e^{(i+1)}=w_e^{(i)}-s\kappa_e^{(i)}w_e^{(i)}+sw_e^{(i)}\sum_{\tau\in E}\kappa_\tau^{(i)}w_\tau^{(i)},
\end{equation}
where all the symbols have the same meaning as those of (\ref{Ni00}) except for Lin-Lu-Yau's Ricci curvature $\kappa_e^{(i)}$ \cite{Lin-Lu-Yau}.
Clearly, the continuous version of (\ref{Lai-1-0}) is
\begin{equation}\label{Lai-1-1}
\frac{d}{dt}w_e(t)=-\kappa_e(t)w_e(t) +w_e(t)\sum_{\tau\in E}\kappa_\tau w_\tau.
\end{equation}
{It's important to note that in most practical networks-especially those with properly scaled edge weights-the total volume does not approach zero. Regularization serves more as a preventive measure than a practical remedy. The value of regularization lies in its theoretical completeness.} Both approaches of \cite{Ni-Lin,Lai X} have successfully detected communities for various complex networks including Karate Club graph, American football games, Facebook ego network, etc.

Though \cite{Ni-Lin} and \cite{Lai X} are greatly successful in applications, there is no theoretical support. In order to supplement this aspect, very recently, Bai et al \cite{Bai-Lin} considered long time existence and convergence of Ricci flow (\ref{flow}) and normalized Ricci flow (\ref{Lai-1-1}) on weighted graphs $(V,E,w(t))$. According to \cite{Bai-Lin}, Lin-Lu-Yau's Ricci curvature \cite{Lin-Lu-Yau} is an appropriate choice for us when considering the intrinsic metric or curvature on graph. Moreover, they proved that under an exit condition (once there exist some $e\in E$ and some $t\in[0,+\infty)$ satisfying
\begin{equation}\label{exitcondition}
w_e(t)>\rho_e(t),
\end{equation}
the edge $e$ is removed),
the normalized Ricci flow (\ref{Lai-1-1}) has a unique global solution on $t\in[0,+\infty)$. Similar as above, if $e=xy$, then $\rho_e(t)$, the distance between $x$ and $y$ at the time $t$, is defined by
$$\rho_e(t)=\inf_{\gamma}\sum_{\tau\in \gamma}w_\tau(t),$$
where $\gamma$ is taken over all paths connecting $x$ and $y$. For the Ricci flow (\ref{flow}) with Ollivier's Ricci curvature replaced by Lin-Lu-Yau's Ricci curvature, the exit condition can still ensure the global existence and uniqueness of the solution. As an example, if a graph is a finite {\it path}, then the exit condition holds, and thus both the Ricci flow and the normalized Ricci flow have unique global solutions. However, for a general connected weighted graph, the exit condition would not hold. This is exactly one reason why Lai et al \cite{Lai X} did surgery at each iteration. \\

The purpose of this paper is to modify systems (\ref{flow}) and (\ref{Lai-1-1}) so that solutions exist for any initial weight and  all time $t\in[0,+\infty)$.
Our first result reads
\begin{theorem}\label{existence}
Let $(V,E,{\bf w})$ be a connected weighted graph, where $E=\{e_1,e_2,\cdots,e_m\}$ denotes the set of  all edges,
${\bf w}=(w_{e_1},w_{e_2},\cdots,w_{e_m})$, $w_{e_i}$ is the weight of the edge $e_i$, $i=1,2,\cdots,m$.
We denote for each $e_i=x_iy_i$, the distance between $x_i$ and $y_i$ by
$$\rho_{e_i}=\min_{\gamma\in\Gamma_{e_i}}\sum_{e_k\in \gamma}w_{e_k},$$
where $\Gamma_{e_i}$ is a set of all paths connecting $x_i$ and $y_i$, $i=1,2,\cdots,m$.
Then for any initial weight ${\bf{w}}_0=(w_{0,1},w_{0,2},\cdots,w_{0,m})\in\mathbb{R}^m_+$, the Ricci flow
\begin{equation}\label{ric-flow}
\left\{\begin{array}{lll}\frac{d}{dt}w_{e_i}(t)=-\kappa_e(t)\rho_{e_i}(t)\\[1.5ex]
w_{e_i}(0)=w_{0,i}\\[1.5ex]
i=1,2,\cdots,m
\end{array}
\right.
\end{equation}
has a unique solution ${\bf{w}}(t)=(w_{e_1}(t),w_{e_2}(t),\cdots,w_{e_m}(t))$ for all $t\in[0,+\infty)$.
\end{theorem}

 Consider a quasi-normalized Ricci flow
\begin{equation}\label{normlize-2}
\left\{\begin{array}{lll}\frac{d}{dt}w_{e_i}(t)=-\kappa_{e_i}(t)\rho_{e_i}(t)+\frac{\sum_{\tau\in E}\kappa_\tau(t)\rho_\tau(t)}{\sum_{h\in E} w_h}\rho_{e_i}(t)
\\[1.2ex]
w_{e_i}(0)=w_{0,i},\,\,i=1,2,\cdots,m.
\end{array}
\right.
\end{equation}
If $\sum_{h\in E}w_h$, the denominator of the second term on the right side of (\ref{normlize-2}), is replaced by $\sum_{h\in E}\rho_h$, then (\ref{normlize-2}) is a normalization of (\ref{ric-flow}). In this sense, we call (\ref{normlize-2}) a quasi-normalization of (\ref{ric-flow}).
Regarding (\ref{normlize-2}), we have the following:
\begin{theorem}\label{normal-theorem}
Let $(V,E,{\bf w})$ be a connected weighted graph, where $E=\{e_1,e_2,\cdots,e_m\}$ denotes the set of  all edges,
${\bf w}=(w_{e_1},w_{e_2},\cdots,w_{e_m})$, $w_{e_i}$ is the weight of the edge $e_i$, $i=1,2,\cdots,m$. Then
for any initial weight ${\bf{w}}_0=(w_{0,1},w_{0,2},\cdots,w_{0,m})\in\mathbb{R}^m_+$, the quasi-normalized Ricci flow (\ref{normlize-2})
has a unique solution ${\bf{w}}(t)=(w_{e_1}(t),w_{e_2}(t),\cdots,w_{e_m}(t))$ for all $t\in[0,+\infty)$.
\end{theorem}

Both the proof of Theorem \ref{existence} and Theorem \ref{normal-theorem} rely on a key observation, namely
\begin{equation}\label{estimate}
-2\max_{\tau\in E}w_\tau(t)\leq\kappa_e(t)\rho_e(t)\leq 2w_e(t)
\end{equation}
for all $e\in E$ and all $t>0$, which will be seen in Section \ref{sec-3} below. Compared with Bai et al's result \cite{Bai-Lin}, we removed the exit condition (\ref{exitcondition}).

Let $T_{\rm max}(\mathbf{w}_0)$ be the maximum time such that (\ref{flow}) has a
solution on the time interval $[0,T_{\rm max}(\mathbf{w}_0))$. Clearly  (\ref{estimate}) does not exclude the possibility that
up to a subsequence,
\begin{equation}\label{estimate-2}
\kappa_e(t)w_e(t)\rightarrow \infty\quad{\rm as}\quad t\rightarrow T_{\rm max}(\mathbf{w}_0).
\end{equation}
Note that the possibility of (\ref{estimate-2}) is the other reason
why Lai et al \cite{Lai X} did surgery at each iteration. Actually, in \cite{Lai X}, approximately $5$ percent of edges are removed in each iteration. In theory, this could lead to inaccurate community detection in some cases.

{ A discrete version of (\ref{ric-flow}) is
\begin{equation}\label{discrete-0}
\left\{\begin{array}{lll}w_{e_i}(t_{n+1})=w_{e_i}(t_n)-s\kappa_{e_i}(t_n)\rho_{e_i}(t_n)\\[1.5ex]
w_{e_i}(0)=w_{0,i},\,\,\,i=1,2,\cdots,m\\[1.5ex]
t_{n}=ns,\,\,\, 0<s<1/2,\,n=0,1,2,\cdots.
\end{array}
\right.
\end{equation}
Denote $a_n=\min_{1\leq i\leq m}w_{e_i}(t_n)$ and $b_n=\max_{1\leq i\leq m}w_{e_i}(t_n)$. By using (\cite{Bai-Lin},
Lemma 1),
we have $-2w_{e_i}(t_n)\leq -\kappa_{e_i}(t_n)\rho_{e_i}(t_n)\leq 2b_n$ and $(1-2s)a_n\leq (1-2s)w_{e_i}(t_n)\leq w_{e_i}(t_{n+1})\leq (1+2s)b_n$  for all $1\leq i\leq m$. This particularly implies
$$(1-2s)a_n\leq a_{n+1}\leq b_{n+1}\leq (1+2s)b_n,\quad\forall n=0,1,2,\cdots.$$
Then an induction argument leads to $a_0(1-2s)^{n}\leq a_n\leq b_n\leq b_0(1+2s)^n$, and thus
\begin{equation}\label{glob-1}a_0(1-2s)^{n}\leq w_{e_i}(t_n)\leq b_0(1+2s)^n,\quad\forall 1\leq i\leq m, \quad \forall n=0,1,2,\cdots.\end{equation}
Observing that a discrete version of (\ref{normlize-2}) reads
\begin{equation}\label{discrete-2}
\left\{\begin{array}{lll}w_{e_i}(t_{n+1})=w_{e_i}(t_n)+s\left\{-\kappa_{e_i}(t_n)\rho_{e_i}(t_n)
+\frac{\sum_{j=1}^m\kappa_{e_j}(t_n)\rho_{e_j}(t_n)}{\sum_{j=1}^m w_{e_j}(t_n)}\rho_{e_i}(t_n)\right\}\\[1.5ex]
w_{e_i}(0)=w_{0,i},\,\,\,i=1,2,\cdots,m\\[1.5ex]
t_{n}=ns,\,\,\, 0<s<\frac{1}{2(m+1)},\,n=0,1,2,\cdots,
\end{array}
\right.
\end{equation}
we obtain by the same method of deriving (\ref{glob-1}) that
\begin{equation}\label{glob-2}a_0(1-2(m+1)s)^{n}\leq w_{e_i}(t_n)\leq b_0(1+4s)^n,\quad \forall 1\leq i\leq m,\quad \forall n=0,1,2,\cdots.\end{equation}
Therefore, both (\ref{discrete-0}) and (\ref{discrete-2}) have global solutions.
We will apply (\ref{discrete-0}) and (\ref{discrete-2}) to detect communities
for complex networks including Karate Club, American football games, Facebook ego network, as well as artificial networks.
However, in \cite{Bai-Lin,Lai X}, estimates similar to (\ref{glob-1}) and (\ref{glob-2}) are generally no longer valid in situations where ``bad edges" appear.
``Bad edges" are defined as edges \( uv \) that meet either of the following conditions: (I) \( w_{uv}(t) > d(u, v)(t) \) for some \( t \in [0, \infty) \), or (II) \( w_{uv}(t) < mt \) (where \( mt \) denotes a pre-defined merge threshold) for some \( t \in [0, \infty) \).
As a consequence, in our algorithms, we don't need to do surgery at each iteration; instead, we just need to delete some  edges with heavy weights after the final iteration, reducing the complexity of the algorithm. This is different from \cite{Ni-Lin,Lai X}.} Moreover,
 our experimental results
are as impressive as those of \cite{Ni-Lin,Lai X}, especially with better modularity. Algorithm \ref{Algo-1} demonstrates the process of applying (\ref{discrete-0})
to community detection. Here $\kappa_e^i$ denotes Lin-Lu-Yau's Ricci curvature, which is calculated via a star coupling Ricci curvature as Lai et al did in \cite{Lai X}.
\begin{algorithm}[H]\label{Algo-1}
        \caption{Community detection using (\ref{discrete-0})}
        \KwIn{an undirected network \( G = \left( {V,E}\right)  \) ; maximum iteration \( T \) ; step size \( s \) .}
        \KwOut{community detection  results of $G$}
       \For{ \( i = 1,\cdots ,T \)}{
            compute the Ricci curvature \( {\kappa }_{e}^{i} \) ; \\
            \( {w}_{e}^{i + 1} = {w}_{e}^{i} - s \times  \left( {{\kappa }_{e}^{i}\times{\rho}_{e}^{i} }\right)  \)  ;\\
        }
        {\it cutoff}$\leftarrow w_{max}$;\\
        \While{cutoff$> w_{min}$}
        {
        \For{ \( e \in  E \)}{
            \If{$w_e>$ cutoff}{
      remove the edge $e$;
      }
            }
        {\it cutoff}$\leftarrow$  {\it cutoff} $-0.01$;\\
        calculate the Modularity, ARI and NMI of $G$;
        }
\end{algorithm}

The remaining part of this paper will be organized as follows: In Section \ref{S3}, we prove that on a connected weighted graph, both the distance function and the Ricci curvature function are Lipschitz continuous with respect to the weights; In Section
\ref{sec-3}, we prove Theorems \ref{existence} and \ref{normal-theorem}, by combining the results in Section \ref{S3}, several conclusions in \cite{Bai-Lin,Bai-Huang}, the key observation (\ref{estimate}) and some ordinary differential system theory; As applications of (\ref{discrete-0}) and (\ref{discrete-2}), experimental process, results and analysis are provided in Section \ref{sec-4}.

\section{The Lipschitz continuity}\label{S3}
Let $(V,E,\mathbf{w})$ be a weighted graph, where $V$ is the set of all vertices, $E=\{e_1,e_2,\cdots,e_m\}$ is the set of all edges, $e_i=x_iy_i$ is the $i$th edge for each $1\leq i\leq m$, and $\mathbf{w}=(w_{e_1},w_{e_2},\cdots,w_{e_m})$ is the weight on $E$. On this weighted graph, the distance between two vertices $u,v\in V$ is given by
\begin{equation}\label{dis-w}d(u,v)=\inf_{\gamma}\sum_{e\in\gamma} w_{e},\end{equation}
where the infimum takes over all paths $\gamma$ connecting $u$ and $v$.
Let $\kappa(u,v)$ be Lin-Lu-Yau's Ricci curvature \cite{Lin-Lu-Yau}, namely
\begin{equation}\label{LLY-curv}\kappa(u,v)=\lim_{\alpha\rightarrow 1}\frac{1}{1-\alpha}\left[1-\frac{W(\mu_u^\alpha,\mu_v^\alpha)}{d(u,v)}\right],\end{equation}
where $d(u,v)$ is defined as in (\ref{dis-w}), $\alpha\in[0,1]$, $W(\mu_u^\alpha,\mu_v^\alpha)$ denotes the usual transportation  distance between two probability measures
$$\mu_u^\alpha(z)=\left\{\begin{array}{lll}
\alpha&{\rm if}& z=u\\[1.2ex]
(1-\alpha)\frac{w_{zu}}{\sum_{y\sim u}w_{yu}}&{\rm if}& z\sim u\\[1.2ex]
0&{\rm if}& {\rm otherwise},
\end{array}\right.\mu_v^\alpha(z)=\left\{\begin{array}{lll}
\alpha&{\rm if}& z=v\\[1.2ex]
(1-\alpha)\frac{w_{zv}}{\sum_{y\sim v}w_{yv}}&{\rm if}& z\sim v\\[1.2ex]
0&{\rm if}& {\rm otherwise}.
\end{array}\right.$$
To simplify notations, we often write $\rho_e=d(x,y)$ and $\kappa_e=\kappa(x,y)$ if $e=xy\in E$. Set
$$\mathbb{R}^m_+=\left\{{\bf{w}}=(w_1,w_2,\cdots,w_m)\in\mathbb{R}^m: w_i>0,i=1,2,\cdots,m\right\}.$$
Now, for any fixed two vertices $u,v\in V$, we have two functions (still denoted by $d(u,v)$, $\kappa(u,v)$), which map $\mathbf{w}\in\mathbb{R}^m_+$ to $d(u,v)\in[0,+\infty)$ and $\kappa(u,v)\in\mathbb{R}$ respectively.\\

We first have the Lipschitz continuity of $d(u,v)$ with respect to the weight $\mathbf{w}$. In particular, we have the following:

\begin{lemma}\label{distance-1}
For any ${\bf{w}},\widetilde{\bf{w}}\in\mathbb{R}^m_+$, if $d$ and $\widetilde{d}$ are two distance functions
determined by $\bf{w}$ and $\widetilde{\bf{w}}$ respectively, then for any two fixed vertices $u,v\in V$, there holds
$$|d(u,v)-\widetilde{d}(u,v)|\leq \sqrt{m}|{\bf{w}}-\widetilde{\bf{w}}|,$$
where $|{\bf{w}}-\widetilde{\bf{w}}|$ denotes the usual Euclidean norm of ${\bf{w}}-\widetilde{\bf{w}}\in\mathbb{R}^m$.
\end{lemma}

\proof  If $d(u,v)\geq\widetilde{d}(u,v)$, then
by definition of $\widetilde{d}(u,v)$, there exists a path $\gamma_0$ connecting vertices $u$ and $v$ such that
$$\widetilde{d}(u,v)=\sum_{e\in \gamma_0}\widetilde{w}_e,$$
which leads to
\begin{eqnarray*}
0\leq d(u,v)-\widetilde{d}(u,v)\leq\sum_{e\in\gamma_0}w_e-\sum_{e\in \gamma_0}\widetilde{w}_e\leq
\sum_{e\in\gamma_0}|w_e-\widetilde{w}_e|\leq \sqrt{m}|{\bf{w}}-\widetilde{\bf{w}}|,
\end{eqnarray*}
where $m$ is the number of the edges in $E$. If $d(u,v)\leq\widetilde{d}(u,v)$,
there exists a path $\gamma_1$ connecting vertices $u$ and $v$ such that
$${d}(u,v)=\sum_{e\in \gamma_1}{w}_e,$$
which leads to
\begin{eqnarray*}
0\leq \widetilde{d}(u,v)-d(u,v)\leq\sum_{e\in\gamma_1}\widetilde{w}_e-\sum_{e\in \gamma_1}w_e\leq
\sum_{e\in\gamma_1}|w_e-\widetilde{w}_e|\leq \sqrt{m}|{\bf{w}}-\widetilde{\bf{w}}|.
\end{eqnarray*}
Thus we get the desired result. $\hfill\Box$\\

We next prove the Lipschitz continuity of $\kappa_e$ with respect to the weight $\mathbf{w}$ as below. {Note that
our construction of various Lipschitz functions will be based on Lemmas $2$ and $3$ in \cite{Bai-Lin}.}

\begin{lemma}\label{curvature-Lip}
For any edge $e=xy\in E$, as a function of ${{\bf{w}}}=(w_{e_1},w_{e_2},\cdots,w_{e_m})$, the Ricci curvature
$\kappa_e$ (see (\ref{LLY-curv}) above) is locally Lipschitz continuous in $\mathbb{R}^m_+$.
\end{lemma}

\proof  Let ${\bf{w}}=(w_{e_1},\cdots,w_{e_m})$ and $\widetilde{{\bf{w}}}=(\widetilde{w}_{e_1},\cdots,\widetilde{w}_{e_m})$ be two vectors in $\mathbb{R}^m_+$. Fixing the edge $e=xy\in E$, we have that  $\kappa_e$ and $\widetilde{\kappa}_e$
are Ricci curvatures determined by ${\bf{w}}$ and $\widetilde{{\bf{w}}}$ respectively. Clearly, we may assume
\begin{equation}\label{hypo}\Lambda^{-1}\leq w_{e_i}\leq\Lambda,\,\,
\Lambda^{-1}\leq \widetilde{w}_{e_i}\leq\Lambda,\,\,
|w_{e_i}-\widetilde{w}_{e_i}|\leq \delta,\,\,\,i=1,2,\cdots,m, \end{equation}
where $\Lambda$ and $\delta$ are two positive constants.
According to (\cite{Munch}, Theorem 2.1; \cite{Bai-Lin}, Remark 1), one has
$$\kappa_e=\inf_{f\in {\rm Lip}\,1,\,\nabla_{yx}f=1}\nabla_{xy}\Delta f,\quad
\widetilde{\kappa}_e=\inf_{f\in \widetilde{{\rm Lip}}\,1,\,\widetilde{\nabla}_{yx}f=1}\widetilde{\nabla}_{xy}\widetilde{\Delta} f,$$
where
$${\rm Lip}\,1=\left\{f:V\rightarrow \mathbb{R}: |f(u)-f(v)|\leq d(u,v),\forall u,v\in V\right\},$$
$$\widetilde{{\rm Lip}}\,1=\left\{f:V\rightarrow \mathbb{R}: |f(u)-f(v)|\leq \widetilde{d}(u,v),\forall u,v\in V\right\},$$
$$\nabla_{xy}f=\frac{f(x)-f(y)}{d(x,y)},\,\, \widetilde{\nabla}_{yx}f=\frac{f(x)-f(y)}{\widetilde{d}(x,y)},\quad\quad$$
$$\Delta f(x)=\frac{1}{\sum_{u\sim x}w_{xu}}\sum_{v\sim x}w_{xv}(f(v)-f(x))\quad\quad\quad$$
and
$$\widetilde{\Delta} f(x)=\frac{1}{\sum_{u\sim x}\widetilde{w}_{xu}}\sum_{v\sim x}\widetilde{w}_{xv}(f(v)-f(x)).\quad\quad\quad$$

With no loss of generality, we assume $d(x,y)\leq\widetilde{d}(x,y)$. Hereafter we distinguish two cases to proceed.

{\bf Case 1.} $\kappa_e\leq \widetilde{\kappa}_e$.

A direct method of variation implies that ${\kappa}_e$ is achieved by a function $f_1\in {\rm Lip}\,1$. In particular $\nabla_{yx}f_1=1$,
$\kappa_e=\nabla_{xy}\Delta f_1$.
In view of Lemma \ref{distance-1}, we find a constant $C_1$ depending only on $\Lambda$ and $m$ satisfying
$$1-C_1|{\bf{w}}-{\widetilde{\bf{w}}}|\leq \frac{d(u,v)}{\widetilde{d}(u,v)}\leq 1+C_1|{\bf{w}}-{\widetilde{\bf{w}}}|,
\quad\forall u,v\in V.$$
As a consequence
$$|f_1(u)-f_1(v)|\leq (1+C_1|{\bf{w}}-{\widetilde{\bf{w}}}|)\widetilde{d}(u,v), \quad\forall u,v\in V.$$
We may assume $|{\bf{w}}-{\widetilde{\bf{w}}}|<1/C_1$. Define a function $f_1^\ast=(1-C_1|{\bf{w}}-{\widetilde{\bf{w}}}|)f_1$. It then follows
that $f_1^\ast\in \widetilde{{\rm Lip}}\,1$ and
$$f_1^\ast(y)-f_1^\ast(x)\leq (1-C_1^2|{\bf{w}}-{\widetilde{\bf{w}}}|^2)\widetilde{d}(x,y)<\widetilde{d}(x,y).$$
Concerning $\widetilde{\bf{w}}$, in view of (\cite{Bai-Lin}, Lemma 2), we find a function
$\widetilde{f}_1\in \widetilde{{\rm Lip}}\,1$ such that
$$\widetilde{f}_1(y)-\widetilde{f}_1(x)=\widetilde{d}(x,y)$$
and for all $u\in V$,
\begin{equation}\label{difference}|\widetilde{f}_1(u)-{f}_1^\ast(u)|\leq \widetilde{d}(x,y)-(f_1^\ast(y)-f_1^\ast(x)).\end{equation}
In particular, (\ref{difference}) follows from (\cite{Bai-Lin}, page 1734, lines 10-13).

Therefore
$$\widetilde{\kappa}_e=\inf_{f\in \widetilde{{\rm Lip}}\,1,\,\widetilde{\nabla}_{yx}f=1}\widetilde{\nabla}_{xy}\widetilde{\Delta} f
\leq \frac{\widetilde{\Delta}\widetilde{f}_1(x)-\widetilde{\Delta}\widetilde{f}_1(y)}{\widetilde{d}(x,y)},$$
and whence
\begin{eqnarray}\nonumber
\widetilde{\kappa}_e-\kappa_e&=&\inf_{f\in \widetilde{{\rm Lip}}\,1,\,\widetilde{\nabla}_{yx}f=1}\widetilde{\nabla}_{xy}\widetilde{\Delta} f
-\frac{\Delta f_1(x)-\Delta f_1(y)}{d(x,y)}\\\nonumber
&\leq&\frac{\widetilde{\Delta}\widetilde{f}_1(x)-\widetilde{\Delta}\widetilde{f}_1(y)}{\widetilde{d}(x,y)}-\frac{\Delta f_1(x)-\Delta f_1(y)}{d(x,y)}\\
&\leq& \left|\frac{\widetilde{\Delta}\widetilde{f}_1(x)}{\widetilde{d}(x,y)}-
\frac{\Delta f_1(x)}{d(x,y)}\right|+\left|\frac{\widetilde{\Delta}\widetilde{f}_1(y)}{\widetilde{d}(x,y)}-
\frac{\Delta f_1(y)}{d(x,y)}\right|.\label{est-1}
\end{eqnarray}
In view of (\ref{difference}), we have for all $u\in V$,
\begin{eqnarray}\nonumber
  \widetilde{f}_1(u)-\widetilde{f}_1(x)&=&\widetilde{f}_1(u)-f_1^\ast(u)+f_1^\ast(u)-f_1^\ast(x)+f_1^\ast(x)-\widetilde{f}_1(x)\\
  &\leq&f_1^\ast(u)-f_1^\ast(x)+2(\widetilde{d}(x,y)-f_1^\ast(y)+f_1^\ast(x))\label{est-2}
\end{eqnarray}
and
\begin{eqnarray}\label{est-3}
  \widetilde{f}_1(u)-\widetilde{f}_1(x)
  \geq f_1^\ast(u)-f_1^\ast(x)-2(\widetilde{d}(x,y)-f_1^\ast(y)+f_1^\ast(x)).
\end{eqnarray}
It follows from (\ref{est-2}) that
\begin{eqnarray}\nonumber
  \frac{\widetilde{\Delta}\widetilde{f}_1(x)}{\widetilde{d}(x,y)}-\frac{\Delta f_1(x)}{d(x,y)}&=&\sum_{u\sim x}\left(\frac{\widetilde{w}_{xu}(\widetilde{f}_1(u)-\widetilde{f}_1(x))}{\widetilde{d}(x,y)\sum_{z\sim x}\widetilde{w}_{xz}}-\frac{{w}_{xu}({f}_1(u)-{f}_1(x))}{{d}(x,y)\sum_{z\sim x}{w}_{xz}}\right)\\\nonumber
  &\leq&\sum_{u\sim x}\left(\frac{\widetilde{w}_{xu}({f}_1^\ast(u)-{f}_1^\ast(x))}{\widetilde{d}(x,y)\sum_{z\sim x}\widetilde{w}_{xz}}
  -\frac{{w}_{xu}({f}_1(u)-{f}_1(x))}{{d}(x,y)\sum_{z\sim x}{w}_{xz}}\right)\\
  &&+\frac{2}{\widetilde{d}(x,y)}(\widetilde{d}(x,y)-f_1^\ast(y)+f_1^\ast(x)).\label{est-4}
\end{eqnarray}
While (\ref{est-3}) gives
\begin{eqnarray}\nonumber
  \frac{\Delta f_1(x)}{d(x,y)}-\frac{\widetilde{\Delta}\widetilde{f}_1(x)}{\widetilde{d}(x,y)}&\leq&\sum_{u\sim x}\left(\frac{{w}_{xu}({f}_1(u)-{f}_1(x))}{{d}(x,y)\sum_{z\sim x}{w}_{xz}}
  -\frac{\widetilde{w}_{xu}({f}_1^\ast(u)-{f}_1^\ast(x))}{\widetilde{d}(x,y)\sum_{z\sim x}\widetilde{w}_{xz}}\right)\\
  &&+\frac{2}{\widetilde{d}(x,y)}(\widetilde{d}(x,y)-f_1^\ast(y)+f_1^\ast(x)).\label{est-5}
\end{eqnarray}
Combining (\ref{est-4}), (\ref{est-5}), $\nabla_{yx}f_1=1$ and the definition of $f_1^\ast$, we obtain
\begin{eqnarray}\nonumber
\left|\frac{\widetilde{\Delta}\widetilde{f}_1(x)}{\widetilde{d}(x,y)}-\frac{\Delta f_1(x)}{d(x,y)}\right|&\leq&
\sum_{u\sim x}\left|\frac{(1-C_1|{{\bf{w}}}-\widetilde{{\bf{w}}}|)\widetilde{w}_{xu}}{\widetilde{d}(x,y)\sum_{z\sim x}\widetilde{w}_{xz}}-
\frac{{w}_{xu}}{{d}(x,y)\sum_{z\sim x}{w}_{xz}}\right||f_1(u)-f_1(x)|\\\nonumber
&&+\frac{2}{\widetilde{d}(x,y)}\left(\widetilde{d}(x,y)-(1+C_1|{{\bf{w}}}-\widetilde{{\bf{w}}}|)d(x,y)\right).\\\nonumber
&\leq&\sum_{u\sim x}\frac{|\widetilde{w}_{xu}d(x,y)\sum_{z\sim x}w_{xz}-w_{xu}\widetilde{d}(x,y)\sum_{z\sim x}\widetilde{w}_{xz}|}{\widetilde{d}(x,y)d(x,y)(\sum_{z\sim x}\widetilde{w}_{xz})(\sum_{z\sim x}w_{xz})}|f_1(u)-f_1(x)|\\\nonumber
&&+C_1|{{\bf{w}}}-\widetilde{{\bf{w}}}|\frac{|f_1(u)-f_1(x)|}{\widetilde{d}(x,y)}\\\nonumber
&&+\frac{2}{\widetilde{d}(x,y)}\left(\widetilde{d}(x,y)-d(x,y)-C_1|{{\bf{w}}}-\widetilde{{\bf{w}}}|d(x,y)\right).\nonumber
\end{eqnarray}
This together with (\ref{hypo}) and Lemma \ref{distance-1} leads to
\begin{equation}\label{est-7}
\left|\frac{\widetilde{\Delta}\widetilde{f}_1(x)}{\widetilde{d}(x,y)}-\frac{\Delta f_1(x)}{d(x,y)}\right|
\leq C|{{\bf{w}}}-\widetilde{{\bf{w}}}|
\end{equation}
for some constant $C$  depending only on $\Lambda$ and $m$.

Also we derive from (\ref{difference}) that
$$\widetilde{f}_1(u)-\widetilde{f}_1(y)\leq f_1^\ast(u)-f_1^\ast(y)+2(\widetilde{d}(x,y)-f_1^\ast(y)+f_1^\ast(x))$$
and
$$\widetilde{f}_1(u)-\widetilde{f}_1(y)\geq f_1(u)-f_1(y)-2(\widetilde{d}(x,y)-f_1^\ast(y)+f_1^\ast(x))$$
for all $u\in V$. Using the same method as proving (\ref{est-7}), we get
\begin{equation}\label{est-8}
\left|\frac{\widetilde{\Delta}\widetilde{f}_1(y)}{\widetilde{d}(x,y)}-\frac{\Delta f_1(y)}{d(x,y)}\right|
\leq C|{{\bf{w}}-\widetilde{{\bf{w}}}}|,
\end{equation}
where $C$ is a constant depending only on $\Lambda$ and $m$.
Inserting (\ref{est-7}) and (\ref{est-8}) into (\ref{est-1}), we find some constant $C$ depending only on $\Lambda$ and $m$
such that
\begin{equation}\label{est-9}\widetilde{\kappa}_e-\kappa_e\leq C|{{\bf{w}}}-\widetilde{{\bf{w}}}|.\end{equation}

{\bf Case 2.} $\kappa_e>\widetilde{\kappa}_e$.

{The main difference from the previous case is that the construction of Lipschitz functions
here will be based on (\cite{Bai-Lin}, Lemma 3), rather than (\cite{Bai-Lin}, Lemma 2). }

 A direct method of variation implies that $\widetilde{\kappa}_e$
is attained by some function $\widetilde{g}_2\in \widetilde{{\rm Lip}}\,1$ satisfying $\widetilde{\nabla}_{yx}\widetilde{g}_2=1$ and
$\widetilde{\kappa}_e=\widetilde{\nabla}_{xy}\widetilde{\Delta}\widetilde{g}_2$. Set $\widetilde{f}_2(u)=\widetilde{g}_2(u)
-\widetilde{g}_2(x)$ for all $u\in V$. Clearly $\widetilde{f}_2\in \widetilde{{\rm Lip}}\,1$, $\widetilde{\nabla}_{yx}\widetilde{f}_2=1$
and $\widetilde{\kappa}_e=\widetilde{\nabla}_{xy}\widetilde{\Delta}\widetilde{f}_2$.
Hence there exists a constant $C_\ast$ depending only on $\Lambda$ and $m$ such that for all $u\in V$,
\begin{equation}\label{C-bd}
|\widetilde{f}_2(u)|=|\widetilde{f}_2(u)-\widetilde{f}_2(x)|\leq \widetilde{d}(u,x)\leq C_\ast,
\end{equation}
since $\widetilde{f}_2(x)=0$.
Define
\begin{equation}\label{def}f_2^\ast(u)=\left\{\begin{array}{lll}\frac{d(x,y)}{\widetilde{d}(x,y)}\widetilde{f}_2(u)&{\rm if}& u\in \{x,y\}\\
[1.5ex]\frac{1}{1+A|{{\bf{w}}}-\widetilde{{\bf{w}}}|}\widetilde{f}_2(u)&{\rm if}& u\in V\setminus\{x,y\},
\end{array}\right.\end{equation}
where $A$ is a large positive number to be determined later.
Keeping in mind $\widetilde{\nabla}_{yx}\widetilde{f}_2=1$, we notice that
\begin{equation}\label{p-1}
f_2^\ast(y)-f_2^\ast(x)=\frac{d(x,y)}{\widetilde{d}(x,y)}(\widetilde{f}_2(y)-\widetilde{f}_2(x))=d(x,y),
\end{equation}
and that for all $u,v\in V\setminus\{x,y\}$,
\begin{eqnarray}\label{p-2}
|f_2^\ast(u)-f_2^\ast(v)|=\frac{1}{1+A|{{\bf{w}}}-\widetilde{{\bf{w}}}|}|\widetilde{f}_2(u)-\widetilde{f}_2(v)|\leq
\frac{1+C_1|{{\bf{w}}}-\widetilde{{\bf{w}}}|}{1+A|{{\bf{w}}}-\widetilde{{\bf{w}}}|} d(u,v),
\end{eqnarray}
\begin{eqnarray}\label{p-3}
|f_2^\ast(u)-f_2^\ast(x)|=\frac{1}{1+A|{{\bf{w}}}-\widetilde{{\bf{w}}}|}|\widetilde{f}_2(u)-\widetilde{f}_2(x)|\leq
\frac{1+C_2|{{\bf{w}}}-\widetilde{{\bf{w}}}|}{1+A|{{\bf{w}}}-\widetilde{{\bf{w}}}|} d(u,x),
\end{eqnarray}
\begin{eqnarray}\nonumber
|f_2^\ast(u)-f_2^\ast(y)|&=&\left|\frac{\widetilde{f}_2(u)-\widetilde{f}_2(x)}{1+A|{{\bf{w}}}-\widetilde{{\bf{w}}}|}-
\frac{d(x,y)}{\widetilde{d}(x,y)}(\widetilde{f}_2(y)-\widetilde{f}_2(x))\right|\\\nonumber
&\leq&\left|\left(\frac{d(x,y)}{\widetilde{d}(x,y)}-\frac{1}{{1+A|{{\bf{w}}}-\widetilde{{\bf{w}}}}|}\right)\widetilde{f}_2(x)\right|+\\\nonumber
&&\left|\left(\frac{d(x,y)}{\widetilde{d}(x,y)}-
\frac{1}{{1+A|{{\bf{w}}}-\widetilde{{\bf{w}}}}|}\right)\widetilde{f}_2(y)\right|
+\frac{|\widetilde{f}_2(u)-\widetilde{f}_2(y)|}{1+A|{{\bf{w}}}-\widetilde{{\bf{w}}}|}\\\label{p-4}
&\leq& O(|{{\bf{w}}}-\widetilde{{\bf{w}}}|)+\frac{1+C_2|{{\bf{w}}}-\widetilde{{\bf{w}}}|}{1+A|{{\bf{w}}}-\widetilde{{\bf{w}}}|} d(u,y).
\end{eqnarray}
If we first fix $\Lambda>1$, then choose a sufficiently large $A>1$, finally take a sufficiently small $\delta>0$, then we conclude
from (\ref{C-bd})-(\ref{p-4}) that $f_2^\ast\in {\rm Lip}\,1$ and $\nabla_{yx}f_2^\ast=1$. For such $f_2^\ast$,
in view of (\cite{Bai-Lin}, Lemma 3), $\forall a, 0<a<d(x,y)$,
there exists a function $f_2\in {\rm Lip}\,1$ such that $f_2(y)-f_2(x)=d(x,y)$ and
$|{f}_2^\ast(z)-f_2(z)|\leq d(x,y)-a$ for all $z\in V$,
which together with $f_2^\ast(x)=0$ implies for all $u\in V$,
$$f_2^\ast(u)-2(d(x,y)-a)\leq f_2(u)-f_2(x)\leq f_2^\ast(u)+2(d(x,y)-a).$$
Now we choose $a=d(x,y)-|{{\bf{w}}}-\widetilde{{{\bf{w}}}}|$. If we take $\delta<\Lambda^{-1}$, then $0<a<d(x,y)$. As a consequence,
there holds
\begin{equation}\label{est-11}f_2^\ast(u)-2|{{\bf{w}}}-\widetilde{{{\bf{w}}}}|\leq f_2(u)-f_2(x)\leq f_2^\ast(u)+2
|{{\bf{w}}}-\widetilde{{{\bf{w}}}}|.\end{equation}
Note that
\begin{eqnarray}\nonumber
\kappa_e-\widetilde{\kappa}_e&=&\inf_{f\in{\rm Lip}\,1,\nabla_{yx}f=1}\nabla_{xy}\Delta f-\widetilde{\nabla}_{xy}\widetilde{\Delta}\widetilde{f}_2\\\nonumber
&\leq&\nabla_{xy}\Delta f_2-\widetilde{\nabla}_{xy}\widetilde{\Delta}\widetilde{f}_2\\\nonumber
&=&\frac{\Delta f_2(x)-\Delta f_2(y)}{d(x,y)}-\frac{\widetilde{\Delta}\widetilde{f}_2(x)-\widetilde{\Delta}\widetilde{f}_2(y)}
{\widetilde{d}(x,y)}\\
&\leq&\left|\frac{\Delta f_2(x)}{d(x,y)}-\frac{\widetilde{\Delta}\widetilde{f}_2(x)}{\widetilde{d}(x,y)}\right|+
\left|\frac{\Delta f_2(y)}{d(x,y)}-\frac{\widetilde{\Delta}\widetilde{f}_2(y)}{\widetilde{d}(x,y)}\right|.\label{di-1}
\end{eqnarray}
On one hand, we calculate
\begin{eqnarray*}
\frac{\Delta f_2(x)}{d(x,y)}-\frac{\widetilde{\Delta}\widetilde{f}_2(x)}{\widetilde{d}(x,y)}&=&
\sum_{u\sim x}\frac{w_{xu}(f_2(u)-f_2(x))}{d(x,y)\sum_{z\sim x}w_{xz}}-\sum_{u\sim x}\frac{\widetilde{w}_{xu}\widetilde{f}_2(u)}{\widetilde{d}(x,y)
\sum_{z\sim x}\widetilde{w}_{xz}}\\
&=&\sum_{u\sim x}\left(\frac{w_{xu}f_2^\ast(u)}{d(x,y)\sum_{z\sim x}w_{xz}}-
\frac{\widetilde{w}_{xu}\widetilde{f}_2(u)}{\widetilde{d}(x,y)
\sum_{z\sim x}\widetilde{w}_{xz}}\right)+O(|{{\bf{w}}}-\widetilde{{{\bf{w}}}}|)\\
&=&\sum_{u\sim x}\left(\frac{w_{xu}}{d(x,y)\sum_{z\sim x}w_{xz}}-\frac{\widetilde{w}_{xu}}{\widetilde{d}(x,y)
\sum_{z\sim x}\widetilde{w}_{xz}}\right)\widetilde{f}_2(u)\\
&&+\sum_{z\sim x}\frac{w_{xu}(f_2^\ast(u)-\widetilde{f}_2(u))}{d(x,y)\sum_{z\sim x}w_{xz}}+O(|{{\bf{w}}}-\widetilde{{{\bf{w}}}}|)\\
&=& O(|{{\bf{w}}}-\widetilde{{{\bf{w}}}}|),
\end{eqnarray*}
since $(\ref{est-11})$ gives $f_2(u)-f_2(x)=f_2^\ast(u)+O(|{{\bf{w}}}-\widetilde{{{\bf{w}}}}|)$,
(\ref{def}) leads to $f_2^\ast(u)-\widetilde{f}_2(u)=O(|{{\bf{w}}}-\widetilde{{{\bf{w}}}}|)$, (\ref{C-bd}) implies
$\widetilde{f}_2(u)=O(1)$, while (\ref{hypo}) and Lemma \ref{distance-1} yields $\widetilde{d}(x,y)/d(x,y)=1+O(|{{\bf{w}}}-\widetilde{{{\bf{w}}}}|)$.
On the other hand, the last term of the inequality (\ref{di-1}) is also $O(|{{\bf{w}}}-\widetilde{{{\bf{w}}}}|)$ in the same way.
Therefore
$$\kappa_e-\widetilde{\kappa}_e\leq C|{{\bf{w}}}-\widetilde{{{\bf{w}}}}|$$
for some constant $C$ depending only on $\Lambda$, $m$ and $\delta$.

Combining Case 1 and Case 2, we complete the proof of the lemma.$\hfill\Box$

\section{Long time existence}\label{sec-3}
In this section, concerning the global existence and uniqueness of solutions to the modified Ricci flow and quasi-normalized Ricci flow on  an arbitrary weighted graph, we prove Theorems \ref{existence} and \ref{normal-theorem}.\\

{\it Proof of Theorem \ref{existence}.} We divide the proof into two steps.\\

{\bf Step 1.} {\it Short time existence.}

Set $w_i=w_{e_i}$, $i=1,2,\cdots,m$. Given any initial value
${\bf{w}}_0=(w_{0,1},w_{0,2},\cdots,w_{0,m})\in\mathbb{R}^m_+$, the modified Ricci flow
(\ref{ric-flow}) is exactly the ordinary differential system
 \begin{equation}\label{equiv}
 \left\{\begin{array}{lll}
 \frac{d}{dt}{\bf{w}}={\bf f}({\bf{w}})\\[1.5ex]
 {\bf{w}(0)}={\bf{w}_0},
 \end{array}\right.
 \end{equation}
where ${\bf w}=(w_1,w_2,\cdots,w_m)\in\mathbb{R}^m_+$ and
 {\bf{f}} is a vector-valued function defined by
\begin{eqnarray*}
{\bf f}:\mathbb{R}^m_+&\rightarrow& \mathbb{R}^m\\
{\bf{w}}&\mapsto& (\kappa_{e_1}\rho_{e_1},\cdots,\kappa_{e_m}\rho_{e_m}).
\end{eqnarray*}
By Lemmas \ref{distance-1} and \ref{curvature-Lip}, we know that all $\rho_{e_i}$, $\kappa_{e_i}$, $i=1,2,\cdots, m$, are locally Lipschitz continuous in
$\mathbb{R}^m_+$, and so is the vector-valued function ${\bf f}({\bf w})$. According to the ODE theory,
for example (\cite{Wang-Zhou-Zhu-Wang}, Section 6.1.1, page 250), there exists a constant $T_0>0$ such that the system (\ref{equiv}) has a unique
solution ${\bf{w}}(t)\in \mathbb{R}^m_+$ for $t\in[0,T_0]$. \\

{\bf Step 2.} {\it Long time existence.}

Define
$$T^\ast=\sup\{T: (\ref{equiv})\, {\rm has\,a\,unique\,solution\,on\,}[0,T]\}.$$
If $T^\ast<+\infty$, then (\ref{equiv}) has a unique solution ${\bf w}(t)=(w_{e_1}(t),w_{e_2}(t),\cdots,w_{e_m}(t))$ on the time interval
$[0,T^\ast)$; moreover, according to an extension theorem of ordinary differential system (\cite{Wang-Zhou-Zhu-Wang}, Section 6.1.1, page 250), there must hold $\mathbf{w}(t)\rightarrow\partial\mathbb{R}^m_+$ or $|\mathbf{w}(t)|\rightarrow+\infty$ as
$t\rightarrow T^\ast$. As a consequence, we have either
\begin{equation}\label{tend-1}\liminf_{t\rightarrow T^\ast}\min\{w_{e_1}(t),w_{e_2}(t),\cdots,w_{e_m}(t)\}= 0,\end{equation}
or
\begin{equation}\label{tend-2}\limsup_{t\rightarrow T^\ast}\max\{w_{e_1}(t),w_{e_2}(t),\cdots,w_{e_m}(t)\}= +\infty.\end{equation}
Denote $\phi(t)=\min\{w_{e_1}(t),w_{e_2}(t),\cdots,w_{e_m}(t)\}$ and $\Phi(t)=\max\{w_{e_1}(t),w_{e_2}(t),\cdots,w_{e_m}(t)\}$.
Thanks to (\cite{Bai-Lin}, Lemma 1) and \cite{Bai-Huang}, there holds for all $e\in E$,
$$
-\frac{2\Phi(t)}{\rho_e(t)}\leq \kappa_e(t)\leq 2,
$$
and thus
\begin{equation}\label{bound}
-2\rho_e(t)\leq -\kappa_e(t)\rho_e(t)\leq 2\Phi(t).
\end{equation}
The power of (\ref{bound}) is evident. On one hand, since $\rho_e(t)\leq w_e(t)$, the left side of (\ref{bound}) gives
$$w_e^\prime(t)=-\kappa_e(t)\rho_e(t)\geq -2\rho_e(t)\geq -2w_e(t),$$
which implies $w_e(t)\geq w_e(0)e^{-2t}$ for all $t\in[0,T^\ast)$. Hence $\phi(t)\geq \phi(0)e^{-2T^\ast}$ for all
$t\in[0,T^\ast)$, contradicting (\ref{tend-1}). On the other hand, the right side of (\ref{bound}) leads to
$$\frac{d}{dt}\sum_{e\in E}w_e(t)\leq 2m\Phi(t)\leq 2m\sum_{e\in E}w_e(t),\quad\forall t\in[0,T^\ast].$$
It follows that
$$\sum_{e\in E}w_e(t)\leq\left(\sum_{e\in E}w_e(0)\right)e^{2mT^\ast},\quad\forall t\in[0,T^\ast],$$
which together with $\Phi(t)\leq \sum_{e\in E}w_e(t)$ contradicts (\ref{tend-2}). Therefore $T^\ast=+\infty$ and
the proof of the theorem is completed. $\hfill\Box$\\

{\it Proof of Theorem \ref{normal-theorem}.}  Suppose that the Ricci curvature
$\kappa_{e_i}$ and the distance  $\rho_{e_i}$
are determined by the weight ${\bf{w}}=(w_{e_1},w_{e_2},\cdots,w_{e_m})$. Let
${\bf{g}}:\mathbb{R}^m_+\rightarrow \mathbb{R}^m$ be a vector-valued function written as ${\bf{g}}({\bf{w}})=(g_1({\bf{w}}),g_2({\bf{w}}),
\cdots,g_m({\bf{w}}))$, where
\begin{equation}\label{g-i}{{g}_i}({\bf{w}})=-\kappa_{e_i}\rho_{e_i}+\frac{\sum_{\tau\in E}\kappa_\tau\rho_\tau}{\sum_{h\in E} {w_h}}\rho_{e_i},\quad i=1,2,\cdots,m.\end{equation}
From Lemmas \ref{distance-1} and \ref{curvature-Lip}, it follows that ${\bf{g}}$ is locally Lipschitz continuous in $\mathbb{R}^m_+$, i.e. for any bounded domain $\Omega\subset\mathbb{R}^m_+$ with $\overline{\Omega}\subset\mathbb{R}^m_+$,  there would exist a constant $C$ depending only on
$\Omega$ such that
$$|{\bf{g}}({\bf{w}})-{\bf{g}}({\widetilde{\bf{w}}})|\leq C|{\bf{w}}-{\widetilde{\bf{w}}}|,\quad\forall {\bf{w}}, {\widetilde{\bf{w}}}\in\Omega.$$
Hence the existence and uniqueness theorem (\cite{Wang-Zhou-Zhu-Wang}, Section 6.1.1, page 250) gives some $T_1>0$ such that the flow (\ref{normlize-2}) has a unique solution ${\bf{w}}(t)$ on the time interval $[0,T_1]$. Set
$$T^\ast=\sup\{T: (\ref{normlize-2})\, {\rm has\,a\,unique\,solution\,on\,}[0,T]\}.$$
If $T^\ast<+\infty$, then (\ref{normlize-2}) has a unique solution ${\bf w}(t)=(w_{e_1}(t),w_{e_2}(t),\cdots,w_{e_m}(t))$ for all
$t\in[0,T^\ast)$; meanwhile, an extension theorem (\cite{Wang-Zhou-Zhu-Wang}, Section 6.1.1, page 250) leads to
\begin{equation}\label{tend-3}\liminf_{t\rightarrow T^\ast}\min\{w_{e_1}(t),w_{e_2}(t),\cdots,w_{e_m}(t)\}= 0\end{equation}
or
\begin{equation}\label{tend-4}\limsup_{t\rightarrow T^\ast}\max\{w_{e_1}(t),w_{e_2}(t),\cdots,w_{e_m}(t)\}= +\infty.\end{equation}
Write $\phi(t)=\min\{w_{e_1}(t),w_{e_2}(t),\cdots,w_{e_m}(t)\}$ and $\Phi(t)=\max\{w_{e_1}(t),w_{e_2}(t),\cdots,w_{e_m}(t)\}$.
We know from (\cite{Bai-Lin}, Lemma 1) and \cite{Bai-Huang} that for all $1\leq i\leq m$,
$$
-\frac{2\Phi(t)}{\rho_{e_i}(t)}\leq \kappa_{e_i}(t)\leq 2.
$$
Since $0<\rho_{e_i}(t)\leq w_{e_i}(t)$ and $\Phi(t)\leq\sum_{i=1}^mw_{e_i}(t)$ for all $t\in[0,T^\ast)$, it follows that
\begin{equation}\label{est-10}-2\rho_{e_i}(t)\leq -\kappa_{e_i}(t)\rho_{e_i}(t)\leq 2\Phi(t)\end{equation}
and
\begin{equation}\label{est-20}-2m\rho_{e_i}\leq \frac{\sum_{\tau\in E}\kappa_\tau\rho_\tau}{\sum_{h\in E} w_h}\rho_{e_i}\leq 2\rho_{e_i},\quad \forall
1\leq i\leq m.\end{equation}
Noticing that $\mathbf{w}$ is a solution of (\ref{normlize-2}), we have by (\ref{est-10}) and (\ref{est-20}),
$$\frac{d}{dt}\sum_{i=1}^mw_{e_i}=-\sum_{i=1}^m\kappa_{e_i}\rho_{e_i}+\frac{\sum_\tau \kappa_\tau \rho_\tau}{\sum_{h}w_{h}}\sum_{i=1}^m\rho_{e_i}\leq 2(m+1)\sum_{i=1}^mw_{e_i}.$$
This leads to
$$
\Phi(t)\leq \sum_{i=1}^mw_{e_i}(t)\leq e^{2(m+1)T^\ast}\sum_{i=1}^mw_{0,i}
$$
for all $t\in[0,T^\ast)$. This concludes (\ref{tend-4}) can not be true.

On the other hand, it follows from (\ref{est-10}) and (\ref{est-20}) that
$$w^\prime_{e_i}=g_i({\bf{w}})\geq-2(m+1)w_{e_i}\quad i=1,2,\cdots,m.$$
Hence
$$w_{e_i}(t)\geq e^{-2(m+1)t}w_{0,i},\quad i=1,2,\cdots,m.$$
As a consequence
$$\phi(t)\geq e^{-2(m+1)T^\ast}\phi(0),$$
which excludes the possibility of (\ref{tend-3}). Therefore $T^\ast=+\infty$, as we desired. $\hfill\Box$\\

{From the proof of Theorem \ref{existence}, we obtain the following interesting corollary:
\begin{corollary}
   Along the modified Ricci flow (\ref{ric-flow}),  the volume growth  satisfies
   $$ e^{-2(m+1)t}\sum_{i=1}^mw_{0,i} \leq \sum_{i=1}^mw_{e_i}(t)\leq e^{2(m+1)t}\sum_{i=1}^mw_{0,i},\quad \forall t\in[0,+\infty).$$
\end{corollary}}

\section{Evaluation}\label{sec-4}
In this section, we apply (\ref{discrete-0}) and (\ref{discrete-2}) to community detection. To simplify notations, in the subsequent experiments, these two Ricci flows are written as Rho and  RhoN respectively. We have listed the algorithm for Rho in the introduction part, and now give the algorithm for RhoN as below.


\begin{algorithm}[H]{\label{Algo-N2}}
        \caption{Community detection using (\ref{discrete-2})}
        \KwIn{an undirected network \( G = \left( {V,E}\right)  \) ; maximum iteration \( T \) ; step size \( s \) .}
        \KwOut{community detection  results of $G$}

       \For{ \( i = 1,\cdots ,T \)}{
            compute the Ricci curvature \( {\kappa }_{e}^{i} \) ; \\
            $sum_r\leftarrow$ $\sum_{e\in E} w_e^i $;\\
            $sum_{\kappa}\leftarrow$ $\sum_{e\in E} \rho_e^i \times {\kappa }_{e}^{i}$;\\
            \( {w}_{e}^{i + 1} = {w}_{e}^{i} - s \times  \left( {{\kappa }_{e}^{i}\times{\rho}_{e}^{i} }+\rho_e^i\times sum_{\kappa} / sum_r \right)  \)  ;\\

        }
         ${\it cutoff}$ $\leftarrow w_{max}$;\\
        \While{cutoff $> w_{min}$}
        {
        \For{ \( e \in  E \)}{
            \If{$w_e>{\it cutoff}$}{
      remove the edge $e$;
      }
            }
        {\it cutoff}$\leftarrow$  {\it cutoff} $-0.01$;\\
        calculate the Modularity, ARI and NMI of $G$;
        }
\end{algorithm}
\noindent According to  \cite{Bai-Huang,Lai X}, in Algorithms $1$ and $2$, the Ricci curvature is calculated by
$$\kappa_e=\frac{1}{\rho_e}\sup_{B}\sum_{u,v\in V}B(u,v)d(u,v),$$
where $B$ is taken over all star couplings defined as follows. Let $e=xy$, $\mu_x^0$ and $\mu_y^0$ be two probability measures as in Section \ref{S3}.
A coupling $B$ between $\mu_x^0$ and $\mu_y^0$ is a star coupling if the following four properties are satisfied:\\
$(i)$ $B(x,y)>0$, $B(u,v)\leq 0$ for all $u\not=x$ or $v\not=y$;\\
$(ii)$ $\sum_{u,v\in V}B(u,v)=0$;\\
$(iii)$ $\sum_{v\in V}B(u,v)=-\mu_x^0(u)$ for all $u\not=x$;\\
$(iv)$
$\sum_{u\in V}B(u,v)=-\mu_y^0(v)$ for all $v\not=y$.

The main complexity of our algorithms comes from finding the shortest path in the graph and solving a linear programming problem. The run times for these tasks are \(O(|E|\log|V|)\) and \(O(|E|D^3)\) respectively, where \(D\) is the average degree of the network, $|E|$ is the number of all edges and $|V|$ is the number of all vertices of the network. Despite the network's sparse connectivity, where \(|D| \ll |E|\), \(O(|E|D^3)\) often exceeds \(O(|E|\log|V|)\) in most cases. Therefore, the computational complexity of our approach is \(O(|E|D^3)\).

Hereafter, we first introduce commonly used criteria for evaluating community detection algorithms. Then, we perform experiments on several real-world networks (such as Karate \cite{Zachary}, Football \cite{Girvan M}, and Facebook \cite{Jure L}) as well as artificial networks to assess the algorithms. To conduct ablation comparison, Rho (\ref{discrete-0}) and RhoN (\ref{discrete-2}) are tested. In the experiments, we compare three traditional machine learning methods that are not based on neural networks. These methods are the Girvan Newman algorithm based on edge betweenness \cite{Girvan M}, the Greedy Modularity algorithm based on modularity maximization \cite{Clauset-Newman-Moore,Reichardt-Bornholdt}, and the Label Propagation algorithm based on stochastic methods \cite{Cordasco-Gargano}. We also use three different models based on Ricci curvature for comparison: unnormalized discrete Ricci flow model with Ollivier Ricci curvature (DORF) \cite{Ni-Lin}, normalized discrete Ricci flow model with Ollivier Ricci curvature (NDORF), and normalized discrete Ricci flow model with star coupling Ricci curvature (NDSRF) \cite{Lai X}.

\subsection{Criteria}
The Adjusted Rand Index (ARI) \cite{Hubert-Arabie} and Normalized Mutual Information (NMI) \cite{Danon-Guilera-Duch} are two commonly used metrics for evaluating the similarity between two partitions. When one of these partitions represents the true labels of network communities, we can assess the detection accuracy of algorithms by comparing the ARI or NMI values between the ground truth and the predicted labels. Unlike ARI and NMI, Modularity \cite{Clauset-Newman-Moore, M.Newman}, denoted as Q, is used in community detection to evaluate partitions by calculating the proportion of edges within communities without requiring the true partition. It assesses a partition based on the principle that edges are dense within communities and sparse between them.

The calculation methods for ARI and NMI each have their strengths. ARI evaluates the consistency between two partitions by counting the number of correctly classified sample pairs, discounting the impact of random partitions. NMI, grounded in information theory, measures the amount of shared information between two partitions to assess their similarity. Each metric has its advantages in practical applications. For example, ARI performs well when dealing with communities of varying sizes, while NMI remains stable across different numbers of communities. Modularity Q, by assessing the proportion of within-community edges, provides a quality measure that is independent of known labels, making it widely applicable in real-world scenarios.

Let \(\{ U_1, U_2, \ldots, U_\mathcal{I} \} \) and \(\{ V_1, V_2, \ldots, V_\mathcal{J} \} \) be two partitions of the set \( S \) of $n$ vertices (nodes). Let \( n_{ij} =|U_i\cap V_j|\) denote the number of vertices in \( U_i\cap V_j \), while \( a_i \) and \( b_j \) represent the numbers of vertices in \( U_i \) and \( V_j \), respectively. Then the explicit expressions of the above mentioned three criteria are listed below.

\begin{itemize}
    \item \textbf{Adjusted Rand Index}
    \[
    \text{ARI} = \frac{\sum_{i=1}^{\mathcal{I}}\sum_{j=1}^{\mathcal{J}} \binom{n_{ij}}{2} - \left( \sum_{i=1}^{\mathcal{I}} \binom{a_i}{2} \sum_{j=1}^{\mathcal{J}} \binom{b_j}{2} \big/ \binom{n}{2} \right)}{\frac{1}{2} \left( \sum_{i=1}^{\mathcal{I}} \binom{a_i}{2} +
    \sum_{j=1}^{\mathcal{J}} \binom{b_j}{2} \right) - \left( \sum_{i=1}^{\mathcal{I}} \binom{a_i}{2} \sum_{j=1}^{\mathcal{J}} \binom{b_j}{2} \big/ \binom{n}{2} \right)},
    \]
    where \( \binom{n}{2} \) is the number of different pairs from \( n \) vertices, while symbols \( \binom{n_{ij}}{2} \), \( \binom{a_i}{2} \) and \( \binom{b_j}{2} \) have the same meaning. ARI ranges from \(-1\) to \(1\), with larger values indicating higher concordance between the two partitions.

    \item \textbf{Normalized Mutual Information}
    \[
    \text{NMI} = \frac{-2 \sum_{i=1}^\mathcal{I} \sum_{j=1}^\mathcal{J} n_{ij} \log \left( \frac{n_{ij} n}{a_i b_j} \right)}{\sum_{i=1}^\mathcal{I} a_i \log \left( \frac{a_i}{n} \right) + \sum_{j=1}^\mathcal{J} b_j \log \left( \frac{b_j}{n} \right)}.
    \]
    NMI ranges from \(0\) to \(1\), with higher values indicating greater similarity between the partitions.

    \item \textbf{Modularity}
    \[
    Q = \sum_{k=1}^N \left( \frac{L_k}{|E|} - \gamma \left( \frac{D_k}{2|E|} \right)^2 \right),
    \]
 where \( N \) represents the number of communities, \( L_k \) is the number of connections within the $k$th community, \( D_k \) is the total degree of vertices in the $k$th community, and \( \gamma \) is a resolution parameter, with a default value of \(1\). The value of \( Q \) ranges from \(-0.5\) to \(1\), with values closer to \(1\) indicating a stronger community structure and better division quality.
\end{itemize}

\subsection{Results}
In this subsection, we outline the model networks and real-world datasets that were used to assess the accuracy of our community detection method.
For the model network, we evaluate a standard and widely used stochastic block model (SBM) \cite{Karrer B} that provides community labels.
As for real-world datasets, we select three different community graphs that come with ground-truth community labels.

\subsubsection{Real world datasets}
Basic information for real-world networks is listed in Table\ref{tab:1}.
\begin{table}[htbp]
\centering
\caption{\label{tab:1}Summary of real-world network characteristics}
\begin{tabular}{ccccccc}
\toprule
networks & vertexes & edges & \#Class & AvgDeg & density  &Diameter\\
\midrule
Karate   & 34      & 78    & 2                     & 4.59   & 0.139  &5 \\
Football & 115     & 613   & 12                    & 10.66  & 0.094   &4\\
Facebook & 775     & 14006 & 18                    & 36.15  & 0.047    &	8\\
\bottomrule
\end{tabular}
\end{table}

The Karate dataset \cite{Zachary} consists of a karate club network with 34 members and 78 edges. The vertices represent members, and the edges represent the connections between members. The actual community structure of the network comprises members from two karate clubs.

For the football dataset \cite{Girvan M}, it is based on the American college football league, which comprises 115 teams (vertices) and 613 matches (edges). The vertices correspond to the teams, while the edges represent the matches between the teams.

The Facebook network \cite{Jure L} is a real dataset from the Stanford Network Analysis Project, representing interaction networks from an online social networking site. In the interaction network, the benchmark community ground truth is organized by well-defined themes such as interests and affiliations.

\begin{figure*}[htbp]
	\centering
         \subfigure[]{
     \includegraphics[scale=0.47]{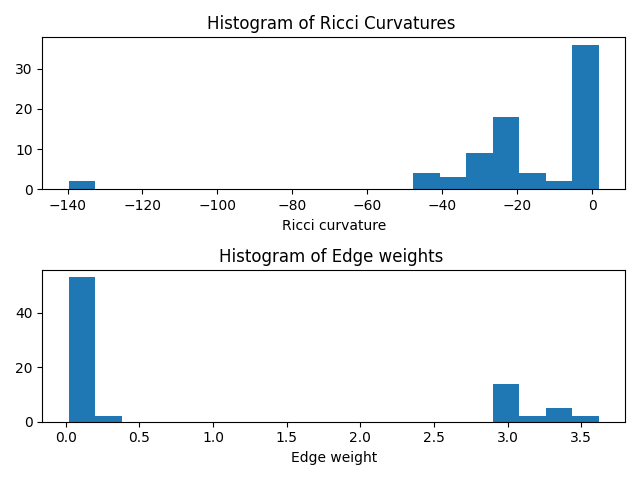}\label{fig:karate_1}
    }
    \subfigure[]{
    \includegraphics[scale=0.47]{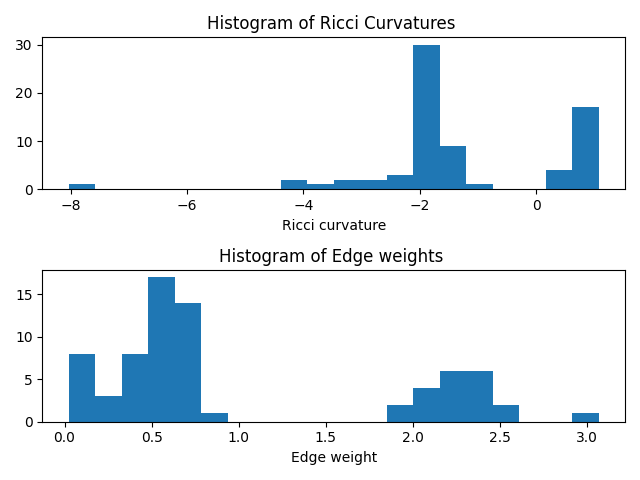}\label{fig:karate_2}
    }
    \subfigure[]{
    \includegraphics[scale=0.6]{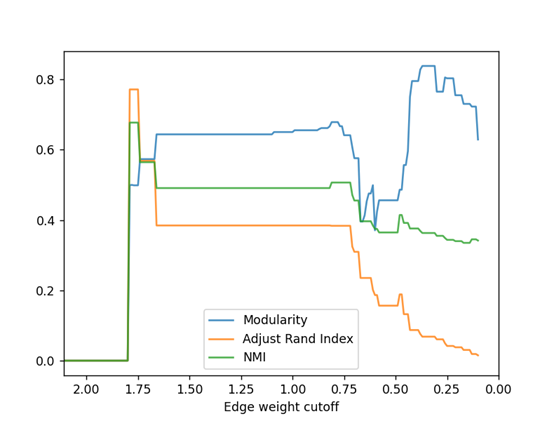}\label{fig:karate_3}
    }	
 	\caption{Karate (a) is the histograms of Ricci Curvatures and Edge Weights before discrete Ricci Flow. (b) is the histograms of Ricci Curvatures and Edge Weights after discrete Ricci Flow. (c) is the evaluation metrics after surgery.}
	\label{fig1}
\end{figure*}

\begin{figure*}[htbp]
	\centering
        \subfigure[]{
     \includegraphics[scale=0.47]{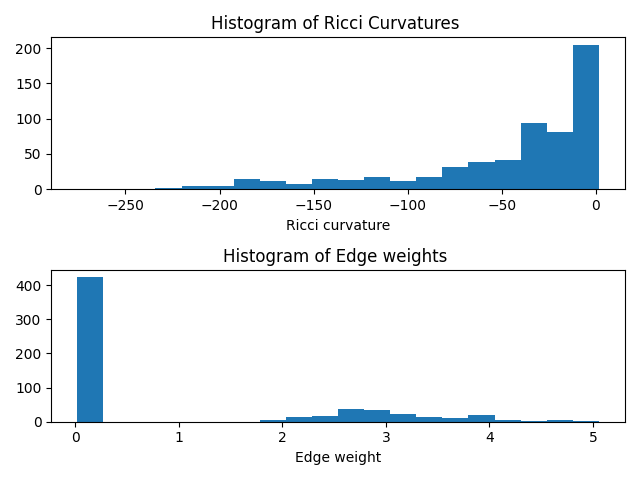}\label{fig:football_1}
    }
    \subfigure[]{
    \includegraphics[scale=0.47]{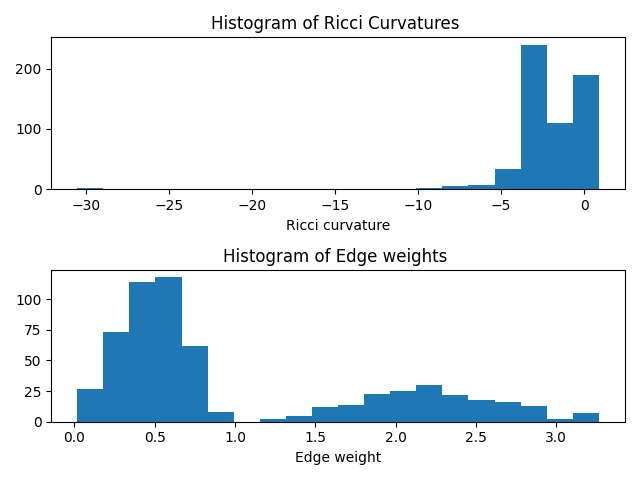}\label{fig:football_2}
    }
    \subfigure[]{
    \includegraphics[scale=0.55]{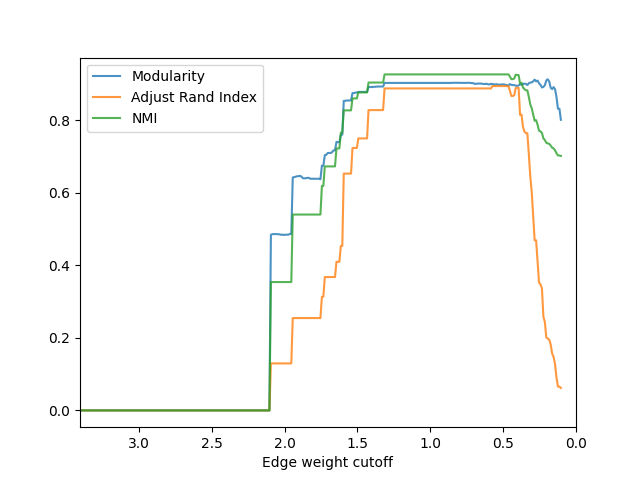}\label{fig:football_3}
    }	
 	\caption{Football (a) is the histograms of Ricci Curvatures and Edge Weights before discrete Ricci Flow. (b) is the histograms of Ricci Curvatures and Edge Weights after discrete Ricci Flow. (c) is the evaluation metrics after surgery.}
	\label{fig2}
\end{figure*}

\begin{figure*}[htbp]
	\centering
        \subfigure[]{
     \includegraphics[scale=0.47]{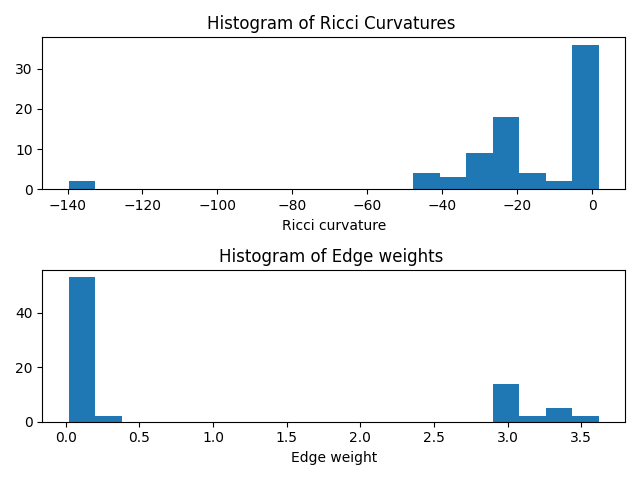}\label{fig:facebook_1}
    }
    \subfigure[]{
    \includegraphics[scale=0.47]{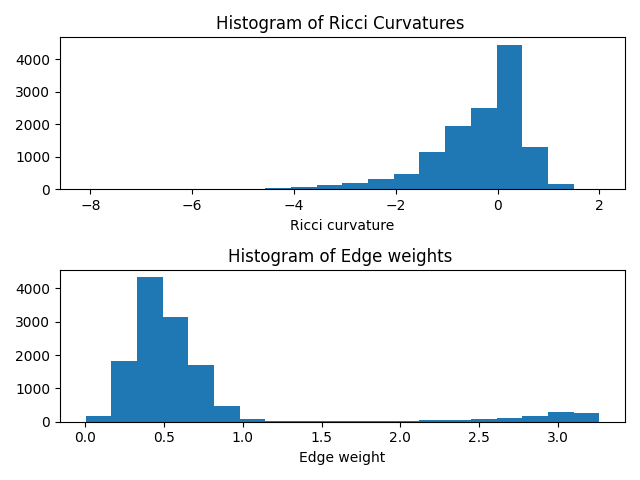}\label{fig:facebook_2}
    }
    \subfigure[]{
    \includegraphics[scale=0.65]{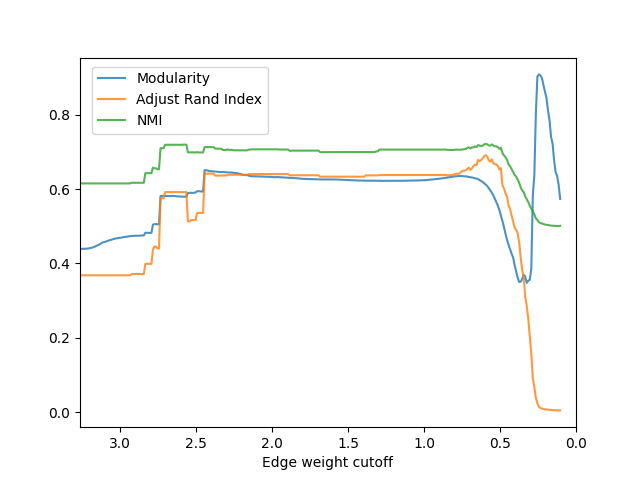}\label{fig:facebook_3}
    }	

 	\caption{Facebook (a) is the histograms of Ricci Curvatures and Edge Weights before discrete Ricci Flow. (b) is the histograms of Ricci Curvatures and Edge Weights after discrete Ricci Flow. (c) is the evaluation metrics after surgery.}
	\label{fig3}
\end{figure*}

The three real-world datasets represent network data at different scales. According to Table \ref{tab:1}, the square of the average degree (denoted as \( D \)) exceeds \( \log \left| V\right|  \). As a result, the computational complexity is \( O\left( {\left| E\right| {D}^{3}}\right) \). The results from our algorithm 2, RhoN (\ref{discrete-2}), for various real-world datasets are presented in sequential order in Figures \ref{fig1}, \ref{fig2}, and \ref{fig3}, corresponding to Karate, Football, and Facebook, respectively. The experimental results of our algorithm 1, Rho (\ref{discrete-0}), and other existing algorithms (Girvan-Newman \cite{Girvan M}, Greedy Modularity \cite{Clauset-Newman-Moore,Reichardt-Bornholdt}, Label Propagation \cite{Cordasco-Gargano}, DORF \cite{Ni-Lin}, NDORF \cite{Bai-Huang} and NDSRF \cite{Lai X}) are presented in Table \ref{tab:2}.

 Since the explanations of Figures 1, 2 and 3 are completely analogous, we  only explain Figure 2 for the Football network, and leave the explanations for Figures 1 and 3 to interested readers.

\begin{itemize}
    \item Figure \ref{fig:football_1} - Before quasi-normalized Ricci flow:

    \textbf{Histogram of Ricci Curvatures:} The distribution is heavily skewed towards negative values, with the bulk of edges showing curvatures from -235 to 0. A pronounced peak near zero indicates minimal contribution from many edges to the total curvature, whereas a tail extending to -235 highlights edges with significantly negative curvature.

    \textbf{Histogram of Edge weights:} There is a notable concentration of edge weights at the lower end, primarily at weight 1, suggesting uniformity across the graph's edges initially.
    \item Figure \ref{fig:football_2} - After quasi-normalized Ricci flow:

    \textbf{Histogram of Ricci Curvatures:} The curvature distribution becomes more balanced and less negatively skewed post-flow, predominantly ranging between -10 and 0, indicative of a normalizing effect of the Ricci flow on the graph's curvatures.

    \textbf{Histogram of Edge weights:} The edge weights now span a broader range from 0 to 3.3, indicating differentiation in edge weights as a result of the Ricci flow process.
    \item Figure \ref{fig:football_3} - Evaluation Metrics Post-Procedure:

    The graph below shows Modularity, Adjusted Rand Index, and NMI plotted against edge weight cutoffs. Initially, all these indicators are zero because only a few edges are removed, and there is minimal community structure. However, the indicators gradually increase, reach a peak, and then stabilize as the {\it cutoff} value increases. When the cutoff approaches $w_{min}$, most of the edges are deleted, leading to a rapid drop in these indicators, essentially reaching zero. This indicates that the Ricci flow has optimized the graph structure up to this point, improving its functionality for community detection.
\end{itemize}

Now we compare the performance of our algorithms with other six algorithms (Girvan-Newman \cite{Girvan M}, Greedy Modularity \cite{Clauset-Newman-Moore,Reichardt-Bornholdt}, Label Propagation \cite{Cordasco-Gargano}, DORF \cite{Ni-Lin}, NDORF \cite{Bai-Huang} and NDSRF \cite{Lai X}) across three real-world network datasets using different metrics. The results in Table \ref{tab:2} demonstrate the effectiveness of our algorithm in detecting communities in real-world scenarios.
\begin{table}[htbp]
\centering
\caption{\label{tab:2}Performance of various algorithms on real-world networks}
\begin{tabular}{cccccccccc}
\toprule
Methods\textbackslash{}Networks & \multicolumn{3}{c}{Karate club} & \multicolumn{3}{c}{Football} & \multicolumn{3}{c}{Facebook} \\
                                                 & ARI       & NMI      & Q        & ARI      & NMI     & Q       & ARI      & NMI     & Q       \\
\midrule
Girvan Newman                                    & \textbf{0.77} & \textbf{0.73} & 0.48          & 0.14          & 0.36          & 0.50 & 0.03          & 0.16          & 0.01          \\
Greedy Modularity                                & 0.57          & 0.56          & 0.58          & 0.47          & 0.70          & 0.82 & 0.49          & 0.68          & 0.55          \\
Label Propagation                                & 0.38          & 0.36          & 0.54          & 0.75          & 0.87          & 0.90 & 0.39          & 0.65          & 0.51          \\
DORF                                             & 0.59          & 0.57          & 0.69          & \textbf{0.93} & \textbf{0.94} & 0.91 & 0.67          & \textbf{0.73} & 0.68          \\
NDORF                                            & 0.59          & 0.57          & 0.69          & \textbf{0.93} & \textbf{0.94} & 0.91 & 0.68          & \textbf{0.73} & 0.68          \\
NDSRF                                            & 0.59          & 0.57          & 0.68          & \textbf{0.93} & \textbf{0.94} & 0.91 & 0.68          & \textbf{0.73} & 0.68          \\
Rho                                     & \textbf{0.77} & 0.68          & 0.82          & 0.89          & 0.92 & 0.90 & 0.64          & 0.72          & 0.63          \\
RhoN                            & \textbf{0.77} & 0.68          & \textbf{0.84} & 0.89          & 0.93 & \textbf{0.92} & \textbf{0.69} & 0.72          & \textbf{0.93}   \\
\bottomrule
\end{tabular}

\end{table}

As shown in Table \ref{tab:2}, in the Karate network, the Rho, RhoN algorithms are slightly less effective than the Girvan Newman algorithm, but only in the NMI metric with a difference of 0.05. When compared to other algorithms, RhoN stands out in the ARI and Q metrics. Notably, RhoN achieves a significantly higher modularity score, averaging 0.3 points higher than the Girvan-Newman, Greedy Modularity, and Label Propagation algorithms. It also surpasses the DORF, NDORF, and NDSRF algorithms by an average of 0.15 points in modularity.

In medium-sized networks such as Football and larger networks like Facebook, the Ricci flow-based algorithms, Rho and RhoN, surpass the three traditional non-neural network-based machine learning methods across all three metrics. In our findings, the $Q$-values are significantly superior, averaging 0.4 points higher than traditional non-neural network-based machine learning methods and 0.25 points higher than Ricci flow-based algorithms. These algorithms effectively capture the community topology of the network and demonstrate superior community detection capabilities regardless of network scale.

If for any three nodes $x$, $y$ and $z$, there always holds $d(x,y)\leq d(x,z)+d(y,z)$, where $d$ is the distance function defined by (\ref{dis-w}), then the network is called satisfying the triangle inequality condition. Compared with other geometry-based discrete Ricci flow algorithms, such as DORF, NDORF, and NDSRF, our algorithms are particularly well-suited for networks without the triangle inequality condition.
\begin{figure}[htbp]
    \centering
    \subfigure{
\begin{minipage}[c]{1\linewidth}
\centering
\includegraphics[width=0.3\linewidth]{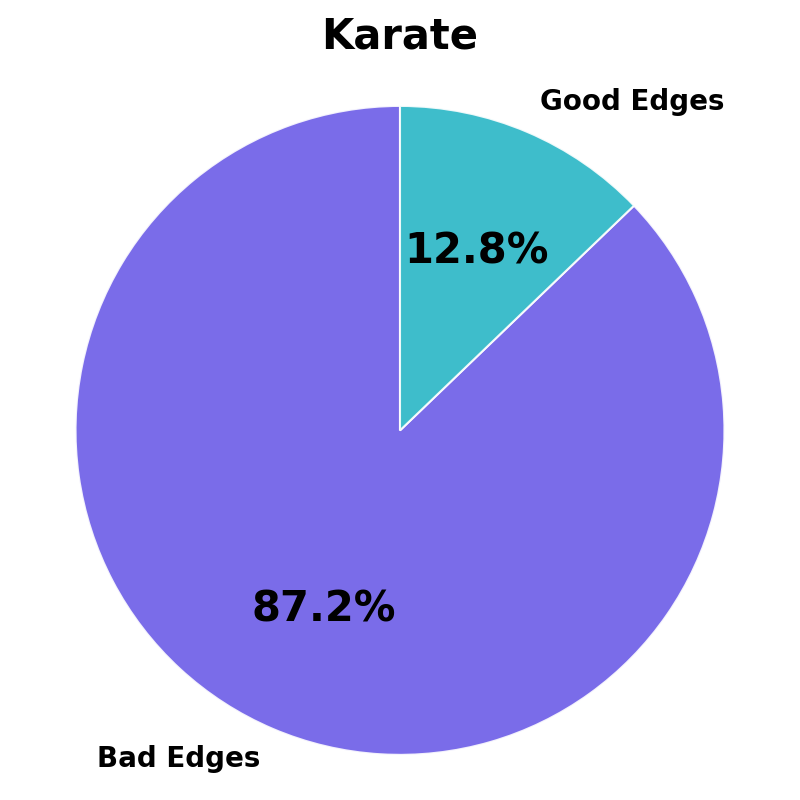}\hspace{4pt}
\includegraphics[width=0.3\linewidth]{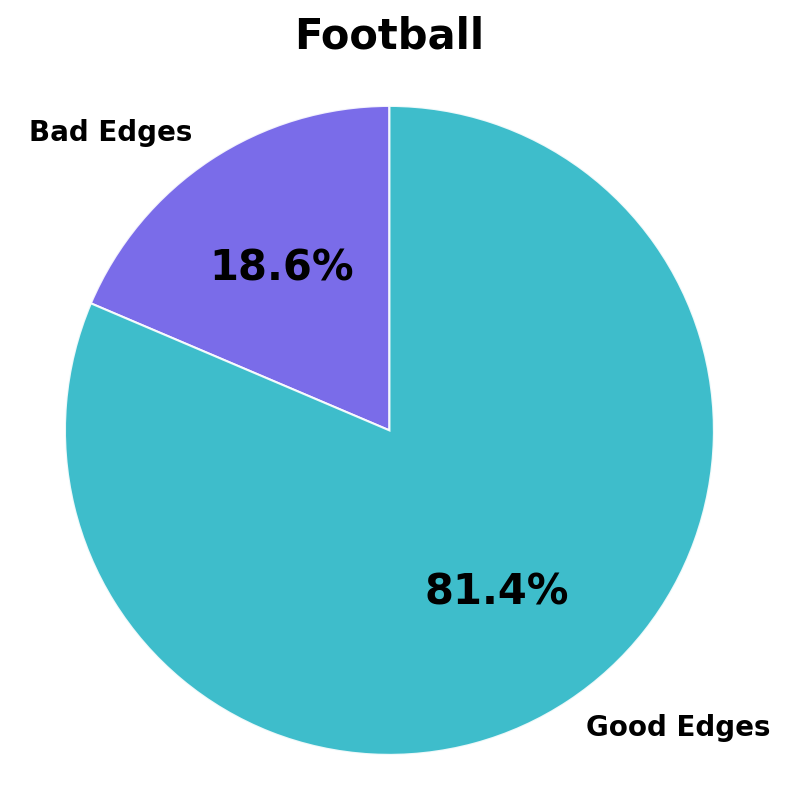}\hspace{4pt}
\includegraphics[width=0.3\linewidth]{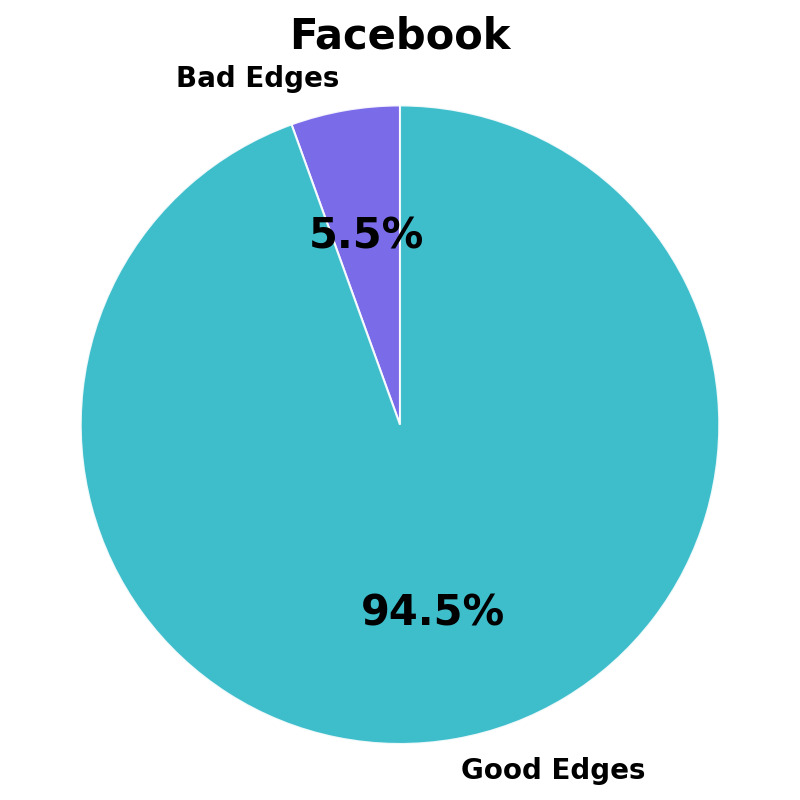}\hspace{4pt}
\end{minipage}
}
    \caption{Frequency of ``bad edges" satisfying the exit condition.}
    \label{fig:frequency}
\end{figure}

{{To provide empirical evidence for the necessity of our proposed modification, we document the frequency of ``bad edges" in real-world networks, as referenced in \cite{Bai-Lin}. Specifically, in the Zachary karate club network, which has 78 edges, 63 of those edges $(80.8\%)$ were classified as bad. In the American college football network, which consists of 613 edges, there were 114 bad edges $(18.6\%)$. Additionally, in the Facebook ego network, which contains 14,006 edges, there were 774 bad edges $(5.5\%)$.

These significant proportions—ranging from mild in large-scale social networks to prevalent in smaller relational graphs—demonstrate that violations of the triangle inequality are neither rare nor negligible in real-world networks. By explicitly addressing these cases, our algorithms fill a critical gap in existing methodologies, enabling a geometric analysis of a broader class of real-world networks. }}

\subsubsection{Artificial datasets}
Stochastic Block Model (SBM) \cite{Karrer B} is a model to generate networks with blocks (communities) randomly. A network generated by SBM is consist with \( n \) nodes, which can be divided into \( k \) blocks, hence there are \( {n}_{i} \) nodes in the \( i \) th block \( {B}_{i} \) . The edges in the graph will be assigned by a probability matrix \( {P}_{k \times  k} = ( {p}_{ij})  \), where \( {p}_{ij} \) presents the edge density from \( {B}_{i} \) to \( {B}_{j} \) . \( {p}_{\text{intra}} \) are usually used as the probability of edges connecting two nodes in the same block, while \( {p}_{\text{inter }} \) for nodes in different blocks. Namely, \( {p}_{\text{intra }} \mathrel{\text{:=}} {p}_{ii},{p}_{\text{inter }} \mathrel{\text{:=}} {p}_{ij}, i \neq  j \), for \( i,j = 1,\ldots ,k \). For networks with implicit community structures, it is necessary that $p_{\text{intra}}$ is greater than $p_{\text{inter}}$, indicating that there are denser connections within the same community and sparser connections between different communities. The upcoming experiments will maintain this implicit community structure.

The synthetic datasets consist of two sets of SBM benchmark networks with different parameters, D1 and D2. The specific parameters and settings are detailed in Table \ref{tab:3}, with any unspecified parameters set to the default values of the SBM tool. To strike a balance between computational speed and intra-community edge density, $p_{\text{intra}}$ is fixed at 0.15. The D1 network is designed to study the effects of varying $p_{\text{inter}}$ values on the accuracy of the algorithm. The mixing parameter $p_{\text{inter}}$ indicates the strength of connections between communities; the larger the value, the denser the connections between nodes of different communities, resulting in more compact communities, a less clear community structure, and a more challenging task for the algorithm to correctly identify the communities.

The D2 network, on the other hand, focuses on studying the effects of different network sizes on the accuracy of the algorithm. A larger size value corresponds to a larger network scale.

\begin{table}[htbp]
\centering
\caption{\label{tab:3}Parameter settings for artificial data}
\begin{tabular}{cccc}
Parameter & Description                                                               & D1        & D2       \\
\toprule
size      & Number of nodes                                                           & 500       & 500-2000 \\
\midrule
k         & \#Communities                                                     & 2         & 2-8      \\
$p_{intra}$  & $p_{ii}$   & 0.15      & 0.15     \\
$p_{inter}$  & $p_{ij}$ & 0.01-0.10 & 0.05     \\
seed      & Random seed                                                               & 0         & 0        \\
directed  & Whether the edges are directed                                            & False     & False    \\
selfloops & Whether self-loops are allowed                                            & False     & False   \\
\bottomrule

\end{tabular}
\end{table}

In a series of synthetic datasets, we demonstrate the scalability and effectiveness of our algorithms based on Rho and RhoN respectively, by comparing its overall performance with the comparison algorithms Girvan Newman \cite{Girvan M}, Greedy Modularity \cite{Clauset-Newman-Moore,Reichardt-Bornholdt}, and Label Propagation \cite{Cordasco-Gargano}, under different metrics.
Since we find that DORF, NDORF, and NDSRF have similar performance to Rho and RhoN, we only compare our algorithms with Girvan Newman, Greedy Modularity and Label Propagation.
The experimental results are shown in the figures below. To minimize the impact of randomness, each data point in the figure represents an average of 10 runs.

\begin{figure}[htbp]
    \centering
    \subfigure[]{
     \includegraphics[scale=0.5]{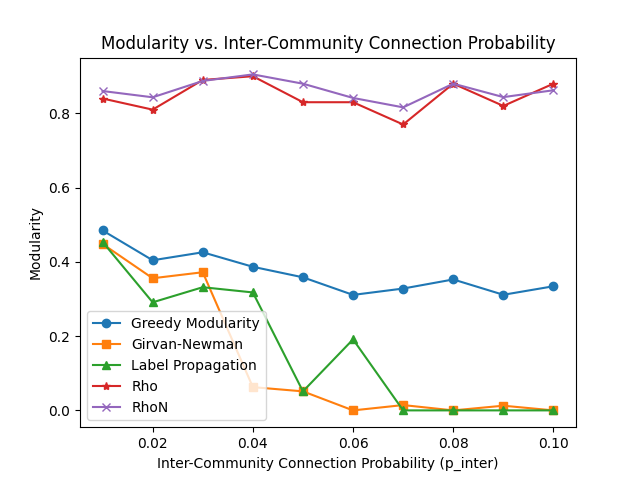}\label{D1_1}
    }
    \subfigure[]{
    \includegraphics[scale=0.5]{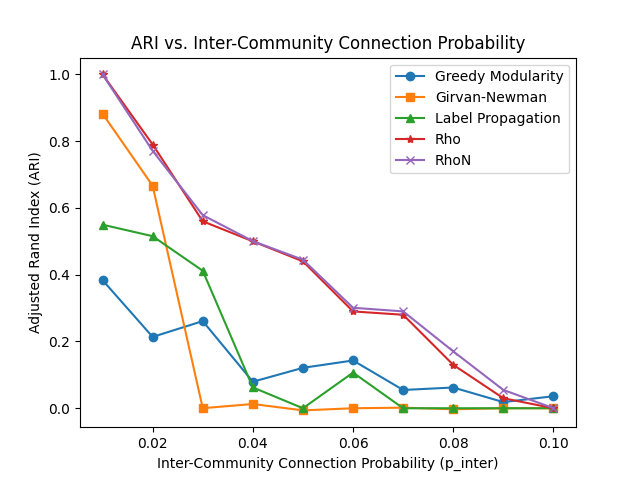}\label{D1_2}
    }
    \subfigure[]{
    \includegraphics[scale=0.5]{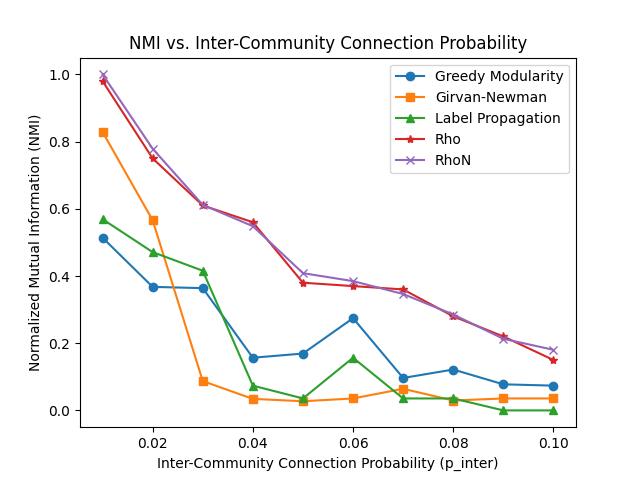}\label{D1_3}
    }

    \caption{Comparing different models on SBM networks (varying $p_{inter}$).}
    \label{fig:D1}
\end{figure}

\begin{figure}[htbp]
    \centering
    \subfigure[]{
     \includegraphics[scale=0.5]{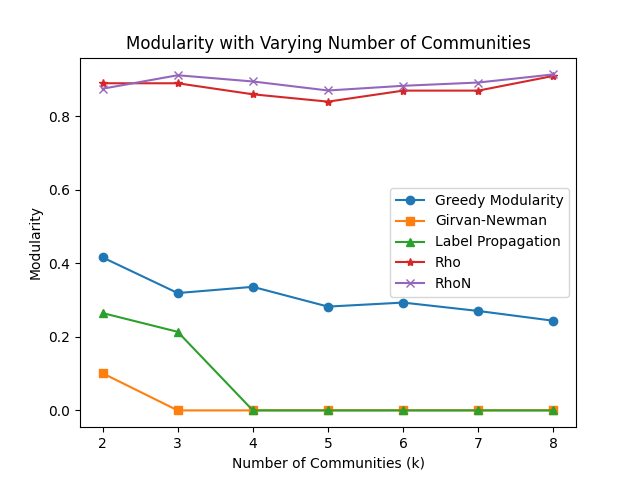}\label{D2_1}
    }
    \subfigure[]{
     \includegraphics[scale=0.5]{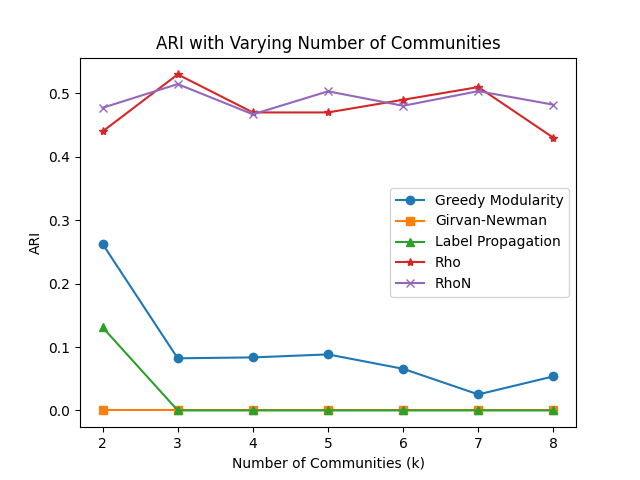}\label{D2_2}
    }
    \subfigure[]{
    \includegraphics[scale=0.5]{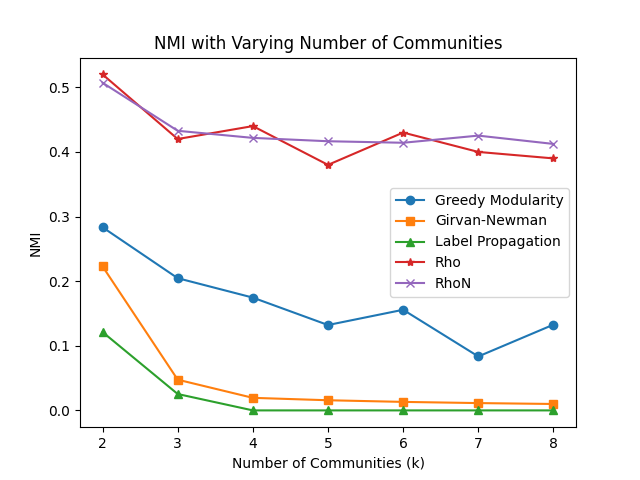}\label{D2_3}
    }

    \caption{Comparing different models on SBM networks (varying $p_{inter}$).}
    \label{fig:D2}
\end{figure}

As the $p_{\text{inter}}$ value increases, the $p_{\text{inter}} / p_{\text{intra}}$ value gradually increases, making the community structure less apparent and the experiments more challenging. In Figure \ref{D1_1}, it is evident that Rho and RhoN show a consistent modularity about 0.8 across all $p_{\text{inter}} $ values, indicating a stable and robust detection of community structures even as inter-community connections increase.  Greedy Modularity also exhibits stability, but with Modularity values around 0.35, which are lower than those of Rho and RhoN. Conversely, Label Propagation and Girvan-Newman show a significant decrease from 0.46 to 0 respectively in Modularity values as the experimental complexity increases.

In Figures \ref{D1_2} and \ref{D1_3}, when $p=0.01$, the experiment presents minimal challenges, and the community structure is quite clear. In this scenario, both Rho and RhoN achieve ARI and NMI values of 1.0, signifying that they accurately identify the network's community structure. As the experimental difficulty increases, all algorithms show a decrease in performance, yet Rho and RhoN consistently maintain the highest ARI and NMI scores across all levels.

In general, Rho and RhoN excel compared to other algorithms across all three metrics, particularly in situations where the community structure is less distinct, indicating superior robustness.

In Figure \ref{D2_1}, Greedy Modularity and Label Propagation both initially have high modularity at k = 2, with values around 0.14 and 0.28 respectively. However, Label Propagation experiences a sharp decline to approximately 0 at k = 3, while Greedy Modularity stabilizes around 0 for larger k. On the other hand, Girvan-Newman starts at 0.4 modularity at k = 2 and maintains this level as k increases.

In Figure \ref{D2_2}, ARI initially starts at approximately 0.27 for $k=2$ but steadily decreases as $k$ increases, eventually nearing 0.05 for higher $k$ values when using the Greedy Modularity algorithm. This trend indicates that Greedy Modularity becomes progressively less effective at accurately identifying community structures as the number of communities grows. The Girvan-Newman algorithm consistently exhibits poor performance, with ARI values remaining close to 0 across all $k$ values, signifying its ineffectiveness in community detection of the Stochastic Block Model. Similarly, Label Propagation demonstrates a sharp decline in ARI after $k=2$, stabilizing near 0, which suggests its inadequacy in handling an increasing number of communities. In contrast, the Rho and RhoN algorithms maintain high and stable ARI values around 0.5, reflecting their robust and consistent ability to detect community structures across varying numbers of communities.

Similarly, Figure \ref{D2_3} reveals that NMI decreases with increasing $k$ values in the Greedy Modularity algorithm, mirroring the trend observed in ARI and highlighting a decline in performance as the number of communities rises. The Girvan-Newman algorithm consistently shows low NMI values across all $k$ levels, with only slightly improving from $k=2$ to $k=3$ before stabilizing at low values. Like its ARI performance, Label Propagation exhibits poor NMI scores, approaching 0 as $k$ increases, further demonstrating its limitations. However, the Rho and RhoN algorithms continue to display relatively high and stable NMI values, ranging from 0.4 to 0.5, underscoring their consistent performance in detecting communities across different $k$ values.

To summarize, though Greedy Modularity and Label Propagation algorithms are effective for networks with a limited number of communities, their performance tends to decline as the number of communities increases. On the other hand, the Girvan-Newman algorithm consistently underperforms across various metrics. In contrast, both algorithms based on Rho and RhoN excel, achieving the highest scores in Modularity, ARI, and NMI, and they also exhibit remarkable stability across networks with different community sizes. This suggests our algorithms are not only robust but also highly capable of accurately detecting community structures in networks with diverse scales and community distributions.

\section{Conclusion}
We have developed two algorithms, which are based on a modified Ricci flow and a quasi-normalized Ricci flow, for weighted graph community detection. Our method significantly improves accuracy and stability compared to other widely used methods, namely, Girvan Newman \cite{Girvan M}, Greedy Modularity \cite{Clauset-Newman-Moore,Reichardt-Bornholdt}, and Label Propagation \cite{Cordasco-Gargano}. It outperforms others in Modularity, ARI, and NMI metrics, showcasing its effectiveness for real-world applications. By cutting edges with large weights only at the final step, our method achieves stable, reliable results, demonstrating robustness, especially in challenging scenarios where community structures are less clear. Future work will aim to optimize these methods further and apply them to larger, more complex networks to enhance performance and applicability in diverse contexts.\\

\noindent
\textbf{Acknowledgements}
We appreciate Lai Xin for her patient and detailed explanation on this topic, especially on her doctoral dissertation
\cite{Lai-doc-thesis}, Liu Shuang for sharing her knowledge about Ricci curvature on connected weighted graph, and referees for their
constructive comments and suggestions on this paper. This research was partly supported by Public Computing Cloud, Renmin University of China.\\

\noindent
\textbf{Data availability}
All data needed are available freely. One can find the codes of our algorithms at https://github.com/mjc191812/Modified-Ricci-Flow.

\section*{Declarations}

\noindent
\textbf{Conflict of interest} The authors declared no potential conflicts of interest with respect to the research, authorship, and publication of this article.\\

\noindent
\textbf{Ethics approval} The research does not involve humans and/or animals. The authors declare that there are no ethics issues to be approved or disclosed.\\



\begin{thebibliography}{00}

\bibitem{Bai-Huang}
S. Bai, A. Huang, L. Lu, S. T. Yau, On the sum of ricci-curvatures for weighted graphs, Pure Appl. Math. Quart. 17 (2021) 1599-1617.

\bibitem{Bai-Lin}
S. Bai, Y. Lin, L. Lu, Z. Wang, S. T. Yau, Ollivier ricci-flow on weighted graphs, Amer. J. Math. 146 (2024) 1723-1747.

\bibitem{Bhowmick}
S. S. Bhowmick, B. S. Seah, Clustering and summarizing protein-protein interaction networks: a survey, IEEE Trans. Knowl. Data Eng. 28 (2015) 638-658.

\bibitem{Brendle-Schoen}
S. Brendle, R. Schoen, Manifolds with \(1/4\)-pinched curvature are space forms, J. Amer. Math. Soc. 22 (2009) 287-307.


\bibitem{Clauset-Newman-Moore}
A. Clauset, M. Newman, C. Moore, Finding community structure in very large networks, Phys. Rev. E 70 (2004) 066111.

\bibitem{Cordasco-Gargano}
G. Cordasco, L. Gargano, Community detection via semi-synchronous label propagation algorithms, BASNA 2010 (2010) 1-8.

\bibitem{Danon-Guilera-Duch}
L. Danon, A. Díaz-Guilera, J. Duch, A. Arenas, Comparing community structure identification, J. Stat. Mech. Theory Exp. 2005 (2005) P09008.

\bibitem{Fortunato}
S. Fortunato, Community detection in graphs, Phys. Rep. 486 (2010) 75-174.

\bibitem{Girvan M}
M. Girvan, M. Newman, Community structure in social and biological networks, Proc. Natl. Acad. Sci. 99 (2002) 7821-7826.

\bibitem{Hamilton}
R. Hamilton, Three-manifolds with positive ricci curvature, J. Differ. Geom. 17 (1982) 255-306.

\bibitem{Hubert-Arabie}
L. Hubert, P. Arabie, Comparing partitions, J. Classif. 2 (1985) 193-218.

\bibitem{Karrer B}
B. Karrer, M. Newman, Stochastic blockmodels and community structure in networks, Phys. Rev. E Stat. Nonlin. Soft Matter Phys. 83 (2011) 016107.

\bibitem{Lai-doc-thesis}
X. Lai, The applications of discrete Ricci curvature
and Ricci flow in graph data analysis, Doctoral dissertation, Renmin University of China, 2023.

\bibitem{Lai X}
X. Lai, S. Bai, Y. Lin, Normalized discrete Ricci flow used in community detection, Phys. A 597 (2022) 127251.

\bibitem{Jure L}
J. Leskovec, SNAP datasets: Stanford large network dataset collection,  http://snap.stanford.edu/data, 2014.

\bibitem{Leskovec}
J. Leskovec, K. Lang, M. Mahoney, Empirical comparison of algorithms for network community detection, Proc. 19th Int. Conf. World Wide Web 2010 (2010) 631-640.

\bibitem{Lin-Lu-Yau}
Y. Lin, L. Lu, S. T. Yau, Ricci curvature of graphs, Tohoku Math. J. 63 (2011) 605-627.

\bibitem{Munch}
F. M\"unch, R. K. Wojciechowski, Ollivier ricci curvature for general graph laplacians: heat equation, laplacian comparison, non-expansion and diameter bounds, Adv. Math. 356 (2019) 11.

\bibitem{Newman M E J}
M. Newman, Modularity and community structure in networks, Proc. Natl. Acad. Sci. 103 (2006) 8577-8582.

\bibitem{M.Newman}
M. Newman, Networks: an introduction, Oxford Univ. Press, 2010.

\bibitem{Ni-Lin}
C. C. Ni, Y. Y. Lin, F. Luo, J. Gao, Community detection on networks with ricci flow, Sci. Rep. 9 (2019) 9984.

\bibitem{Ni C C}
C. C. Ni, Y. Y. Lin, F. Luo, J. Gao, X. Gu, E. Saucan, Ricci curvature of the internet topology, 2015 IEEE Conf. Comput. Commun. INFOCOM (2015) 2758-2766.

\bibitem{Ollivier-1}
Y. Ollivier, Ricci curvature of markov chains on metric spaces, J. Funct. Anal. 256 (2009) 810-864.

\bibitem{Ollivier-2}
Y. Ollivier, Ricci curvature of metric spaces, C. R. Math. 345 (2007) 643-646.

\bibitem{Peel}
L. Peel, D. Larremore, A. Clauset, The ground truth about metadata and community detection in networks, Sci. Adv. 3 (2016) 08.

\bibitem{Perelman}
G. Perelman, The entropy formula for the ricci flow and its geometric applications, arXiv: math/0211159, 2002.

\bibitem{Reichardt-Bornholdt}
J. Reichardt, S. Bornholdt, Statistical mechanics of community detection, Phys. Rev. E 74 (2006) 016110.

\bibitem{Tauro S L}
S. Tauro, C. Palmer, G. Siganos, et al., A simple conceptual model for the internet topology, GLOBECOM'01 IEEE Global Telecommun. Conf. 3 (2001) 1667-1671.

\bibitem{Wang-Zhou-Zhu-Wang}
G. Wang, Z. Zhou, S. Zhu, S. Wang, Ordinary differential equations (in Chinese), Higher Education Press, 2006.

\bibitem{Wasserman S}
S. Wasserman, K. Faust, Social network analysis: methods and applications,
Cambridge University Press, 1994.

\bibitem{Yang-Algesheimer}
Z. Yang, R. Algesheimer, C. Tessone, A comparative analysis of community detection algorithms on artificial networks, Sci. Rep. 6 (2016) 30750.

\bibitem{Zachary}
W. Zachary, An information flow model for conflict and fission in small groups, J. Anthropol. Res. 33 (1977) 452-473.

\end{thebibliography}
\end{document}